\documentclass[11pt,a4paper]{article}
\usepackage{amssymb,amsmath,amsthm,amsfonts,mathrsfs,bm,enumerate,esint}
\allowdisplaybreaks[1]
\numberwithin{equation}{section}

\usepackage[colorlinks=true, pdfstartview=FitV, linkcolor=blue, citecolor=red, urlcolor=blue,pagebackref=false]{hyperref}

\usepackage{tikz}
\usepackage{float}
\usepackage[font=footnotesize]{caption}

\usepackage{color} 
\newtheorem{theorem}{Theorem}[section]

\newtheorem{corollary}[theorem]{Corollary}
\newtheorem{lemma}[theorem]{Lemma}
\newtheorem{proposition}[theorem]{Proposition}

\newtheorem{remark}[theorem]{Remark}

\newcommand{\R}{\mathbb R}

\newcommand{\Z}{\mathbb Z}
\newcommand{\N}{\mathbb N}
\newcommand{\T}{\mathbb T}
\newcommand{\eps}{\varepsilon}
\DeclareMathOperator{\diver}{div}
\DeclareMathOperator{\supp}{supp}

\newcommand{\one}{\mathbf 1}
\newcommand{\X}{\mathcal X}
\newcommand{\U}{\mathcal A}

\newcommand{\cJ}{\mathcal J}

\newcommand{\G}{\mathcal G}
\newcommand{\A}{\mathcal A}
\newcommand{\cH}{\mathcal H}
\newcommand{\cN}{\mathcal N}
\newcommand{\nlam}{\nabla_\lambda}

\newcommand{\dd}{\, \mathrm{d}}

\newcommand{\anorm}[2]{\left\|#1\right\|_{\underline L^2(#2)}}
\newcommand{\pnorm}[3]{\left\|#1\right\|_{\underline L^{#2}(#3)}}
\newcommand{\Dini}{\mathcal D}

\usepackage{comment}

\begin{document}

\title{
Quantitative estimates in quasiperiodic elliptic homogenization with arbitrary frequencies
}

\author{
Weisheng Niu \thanks{School of Mathematical Science, Anhui University,
Hefei 230039,  China. 
\textit{Email address: } \texttt{niuwsh@ahu.edu.cn}}
\and
Yousheng Wang \thanks{School of Mathematical Sciences, University of Chinese Academy of Sciences, Beijing, 100049, China.
\textit{Email address: }\texttt{wangyousheng24@mails.ucas.ac.cn}}
\and
Yao Xu \thanks{School of Mathematical Sciences, University of Chinese Academy of Sciences, Beijing, 100049, China.
\textit{Email address: }\texttt{xuyao@ucas.ac.cn}}
\and
Jinping Zhuge \thanks{Corresponding author. Morningside Center of Mathematics, Academy of Mathematics and Systems Science, Chinese Academy of Sciences, Beijing 100190, China.
\textit{Email address: }\texttt{jpzhuge@amss.ac.cn}} 
}

\date{}

\maketitle

\begin{abstract}
    We develop a frequency decomposition approach to study quantitative estimates for strongly elliptic equations with quasiperiodic coefficients $A(\lambda x)$ and arbitrary frequency matrices $\lambda \in \R^{N\times d}$. We prove uniform Lipschitz estimates with constants independent of the entries of $\lambda$, which particularly include equations of finitely many periodic scales without any scale-separation condition. In addition, we construct correctors and flux correctors with quantitative sublinear growth rates explicitly depending on $\lambda$ and the regularity exponent of the coefficient matrix $A$. Examples under Diophantine conditions show that these growth rates are generally optimal. Both results are established within the same framework of Littlewood-Paley decompositions for lifted correctors and flux correctors on the torus $\T^N$. In addition to the classical Littlewood-Paley decomposition that captures regularity in $\T^N$, we also introduce a directional Littlewood-Paley decomposition which captures the oscillations in physical space $\R^d$ to handle the degenerate quasiperiodic structure of corrector equations. 
    Our quantitative framework is new and robust for quasiperiodic homogenization.
  \end{abstract}

\noindent {\bf Keywords.} quasiperiodic homogenization, Lipschitz estimate, correctors, Littlewood-Paley decomposition

\noindent {\bf MSC2020.} 35B27, 74Q05, 35B65. 

\tableofcontents

\section{Introduction}
\label{sec.intro}

\subsection{Motivations}
This paper is devoted to quantitative estimates of elliptic equations in divergence form with quasiperiodic coefficients. Precisely, letting $d \ge 2$ and $N \ge 1$, we consider  a family of divergence-form elliptic operators
\begin{equation}\label{eq:operator}
\mathcal L_\lambda:=-\diver(A_\lambda\nabla),\qquad \text{in } \R^d,
\end{equation}
where $\lambda \in \R^{N\times d}$ is an arbitrary frequency matrix,
$A_\lambda(x):=A(\lambda x)$ and
 $A:\R^N\to\R^{d\times d}$ is a real, $\Z^N$-periodic coefficient matrix satisfying the strong ellipticity condition. The exact ellipticity and regularity assumptions are stated in Subsection~\ref{subsec:main-results}. The novelty of this paper is that we do not impose any restriction on the frequency matrix $\lambda$ --- it may have any rank, any real-valued entries, with arbitrary sizes spreading across many scales. In particular, this includes multiscale periodic coefficients 
 \begin{equation}\label{eq:multiscale-A}
 A_\lambda(x)
 = A(x/\eps_1, x/\eps_2,\dots, x/\eps_n),
 \end{equation}
 where $\eps_1\geq \cdots\geq \eps_n>0$ are arbitrary scale parameters and
 \begin{equation}\label{eq:def-multiscale}
\lambda=\begin{pmatrix}\varepsilon_1^{-1}I_d\\ \vdots\\ \varepsilon_n^{-1}I_d\end{pmatrix} \in \R^{N \times d}, \qquad N = nd.
\end{equation}
Uniform estimates in $\lambda$ therefore give uniform estimates in $\eps$ without any scale-separation.

For a fixed $\lambda$ and coefficient matrix of the form $A(\lambda x/\eps)$ with one microscopic scale $\eps \ll 1$, the quantitative quasiperiodic or even more general almost-periodic elliptic homogenization goes back to  Kozlov \cite{Koz78}. In the quasiperiodic case, one can lift the corrector equations to the phase torus $\T^N = \R^N/\Z^N$, which makes periodicity available. But the resulting periodic elliptic operator is generally degenerate and depends essentially on the frequency matrix $\lambda$. Therefore, inverting that operator would encounter small divisors. Kozlov obtained bounded quasiperiodic correctors and thus $O(\eps)$ convergence rates under a Diophantine condition on the frequency matrix and a higher-order regularity assumption on the coefficient matrix \cite[Section~3]{Koz78}. For arbitrary $\lambda$, small divisors can instead postpone homogenization to arbitrarily large scales and convergence rates can be arbitrarily slow \cite{BBMM05}. In \cite{S15}, Shen developed a different framework for almost-periodic homogenization without lifting the cell problem to the phase torus. In fact, Shen introduced a modulus $\rho(R)$ (also see earlier analogue in \cite[Section~8]{ACS14} for Hamilton-Jacobi homogenization) to quantify the almost-periodicity of the coefficient matrix $A$ in the physical space $\R^d$, where
\begin{equation}\label{eq:Shen-modulus}
    \rho(R) = \sup_{y\in \R^d} \inf_{z\in B_R} \| A(\cdot + y) - A(\cdot + z) \|_{L^\infty(\R^d)}.
\end{equation}
The decay rate as $R \to \infty$ of this modulus implies suboptimal sublinear growth of correctors and  convergence rates. 
Later, Armstrong-Gloria-Kuusi \cite{AGK16} generalized the above modulus to higher-order $\rho_k(R)$, which can decay as $O(R^{-1-\sigma})$, leading to globally bounded correctors and a sharp convergence rate of $O(\eps)$. Further developments or generalizations under this framework can be found in \cite{AS16,SZ18a,ZG17,XN21,GS25}.

On the other hand, uniform regularity is another central topic in homogenization with the pioneering work of Avallenda-Lin \cite{AL87,AL91} for periodic coefficients with one scale parameter, using a compactness method. In stochastic homogenization, Armstrong-Smart \cite{ASm16} developed a robust quantitative method in the spirit of Campanato iteration to establish the large-scale uniform Lipschitz estimate. Nowadays, this method, further developed in different forms (see e.g. \cite{S17,GNO20}), has become a standard method for uniform Lipschitz estimates in either stochastic or deterministic homogenization. We refer to recent monographs \cite{S18bk,AKMbook1} for details of this method.

For multiscale periodic homogenization in the form of \eqref{eq:multiscale-A}, the special case $\eps_i = \eps^i$ was first studied in \cite{BLP78} and the general qualitative homogenization theorem was proved in \cite{AB96} under a scale-separation condition, namely, $\{\eps_i\}$ and all their successive ratios $\{ \eps_{i+1}/\eps_i \}$ tend to zero. In \cite{NSX20}, Niu-Shen-Xu established quantitative convergence rates, and derived the uniform Lipschitz regularity under an algebraic scale-separation condition $\frac{\eps_{i+1}}{\eps_i} \le (\frac{\eps_i}{\eps_{i-1}})^\sigma$ for some $\sigma>0$. Recently, Niu-Xu-Zhuge \cite{NXZ26} established the optimal convergence rates with a more accurate effective equation for real analytic coefficients and proved the uniform Lipschitz estimate under a mild double-logarithmic scale-separation condition $\frac{\eps_i}{\eps_{i+1}} \gtrsim \log \log \eps_i^{-1}$. Note that these regularity results were obtained under certain scale-separation conditions. This motivates us to ask a natural interesting question: can we show the uniform regularity independent of the scale-separation condition? In \cite{NZ23}, the authors used the compactness method to prove the uniform $C^\alpha$ regularity for any $\alpha \in (0,1)$ and for arbitrary $(\eps_1,\dots, \eps_n) \in (0,1]^n$ without any scale-separation condition, and conjectured that the uniform Lipschitz estimate is true. Later in \cite{NZ24} they developed a quantitative method using Dirichlet's theorem from number theory and a real-variable argument to prove the uniform $W^{1,p}$ estimate for any $p\in (1,\infty)$. This is equivalent to the uniform $W^{1,p}$ estimate for the general quasiperiodic elliptic operator \eqref{eq:operator}; see \cite[Corollary 1.2]{NZ24}.
However, the uniform Lipschitz estimate, as the endpoint situation, is still missing, because Armstrong-Smart approach is very difficult to apply under no scale-separation condition. 

In this paper, we are interested in two fundamental questions in quasiperiodic elliptic homogenization: uniform Lipschitz estimate of solutions and sublinear growth of correctors. We prove the uniform Lipschitz estimate for the elliptic operator \eqref{eq:operator} independent of the entries of $\lambda$, including the special case \eqref{eq:multiscale-A} (they are actually equivalent). This solves the conjecture raised in \cite{NZ23}. Furthermore, for an arbitrary frequency matrix $\lambda$, we will establish the optimal sublinear growth rates of correctors and flux correctors, which would imply the optimal quantitative homogenization at sufficiently large scales. These two results will be established under the same framework using the Littlewood-Paley decomposition and its directional variant applied to the lifted cell problem.

\subsection{Assumptions and main results}
\label{subsec:main-results}

We first state our assumptions on the coefficient matrix $A(Y): \R^N \to \R^{d\times d}$.
\begin{itemize}
    \item Strong ellipticity: there exists $1\le\Lambda<\infty$ such that
    \begin{equation}\label{eq:ellipticity}
        \Lambda^{-1}|\zeta|^2\le\zeta\cdot A(Y)\zeta,\qquad
        |A(Y)\zeta|\le\Lambda|\zeta|\qquad \text{for any } Y \in\R^N,\ \zeta\in\R^d.
\end{equation}
    \item $\Z^N$-periodicity: for all $Z\in \Z^N$ and $Y \in \R^N$,
    \begin{equation}\label{eq:periodicity}
        A(Y+Z)  = A(Y).
    \end{equation}

    \item $C^\alpha$ regularity: there exist $\alpha \in (0,1]$ and $\Lambda_1 > 0$ such that
    \begin{equation}\label{eq:holder}
        |A(Y)-A(Y')|\le \Lambda_1 |Y-Y'|^\alpha \qquad \text{for any } Y,Y'\in\R^N.
    \end{equation}
    \end{itemize}
    We will call $(d, N, \Lambda, \Lambda_1, \alpha)$ the character of $A$.

To state the main result, we also use the following notations for convenience.
For a measurable set $E \subset \R^d$ with $0<|E|<\infty$,  $1\le p<\infty$, set
$$
\fint_E g:=\frac1{|E|}\int_Eg,\qquad
\pnorm{g}{p}{E}:=\left(\fint_E|g|^p\right)^{1/p}.
$$
For $0<\alpha \le1$, let
$$
[g]_{C^{0,\alpha}(E)}:=\sup_{\substack{x,x'\in E \\x \neq x'}}
\frac{|g(x)-g(x')|}{|x-x'|^\alpha}.
$$

For brevity, we write $B_R:=B(x_0, R)$ for the ball of radius $R>0$ centered at some fixed point $x_0\in \mathbb{R}^d$, if $x_0$ makes no difference in the statement. 

\begin{theorem}[Uniform Lipschitz estimate]\label{thm:main}
Let $A$ satisfy \eqref{eq:ellipticity}--\eqref{eq:holder}, $p>d$ and $0<\gamma\le1$. For any $\lambda\in\R^{N\times d}$ and $R>0$, let $u\in H^1(B_R)$ be a weak solution of
\begin{equation}\label{eq:equation-BR}
-\operatorname{div}(A_\lambda\nabla u)=f+\operatorname{div}g\quad\text{in }B_R\,,
\end{equation}
where $f\in L^p(B_R)$ and $g\in C^{0,\gamma}(B_R;\R^{d})$. Then
\begin{equation}\label{eq:Lipschitz}
\|\nabla u\|_{L^\infty(B_{R/2})}\le C\left\{\|\nabla u\|_{\underline L^2(B_R)}+R\|f\|_{\underline L^p(B_R)}+R^\gamma[g]_{C^{0,\gamma}(B_R)}\right\}\,.
\end{equation}
The constant $C$ depends only on $p,\gamma$ and the character of $A$, and is independent of $\lambda,R$.
\end{theorem}

The novelty of the above theorem is the uniformity in arbitrary $\lambda \in \R^{N\times d}$ and $R>0$. This permits nearly resonant quasiperiodic structure with small divisors, and multiscale matrix \eqref{eq:multiscale-A} without any scale-separation condition. Therefore, Theorem \ref{thm:main} solves the conjecture raised in \cite{NZ23} (restated in \cite{NZ24,NXZ26}). We mention that our theorem can be extended to elliptic systems and boundary estimates as well.

To state our second result on growth rates, we first introduce the correctors and flux correctors. The correctors $\chi = (\chi_j)_{1\le j\le d}$ are the solutions of the equations
\begin{equation}\label{eq:corrector}
    -\diver (A_\lambda (e_j + \nabla \chi_j) ) = 0 \qquad \text{ in } \R^d.
\end{equation}
The flux correctors $\phi = (\phi_{kij})_{1\le k,i,j\le d}$ are the solutions of the equations
\begin{equation}\label{eq:fluxcorrector}
    \Delta \phi_{kij}=\partial_k[A_\lambda(I+\nabla\chi)]_{ij}-\partial_i[A_\lambda(I+\nabla\chi)]_{kj}\qquad\text{in }\R^d,
\end{equation}
where $I$ is the identity matrix and $(\nabla \chi)_{ij} = \partial_i \chi_j$. Here we require that $\chi$ and $\phi$ satisfy the sublinear growth at infinity. 

By assuming the correctors have the same quasiperiodic structure as the coefficient matrix $A_\lambda$, we can lift problem \eqref{eq:corrector} to a cell problem on $\T^N$ (see Section \ref{sec.corrector}) with a solution in the form of $u(x) = U(\lambda x)$, where $U(Y)$ is a $\Z^N$-periodic function.
Observe that
\begin{equation}\label{eq:intro-lifting}
    \nabla_x u(x) = \nabla_\lambda U(\lambda x) = \sum_{n\in \Z^N} 2\pi i\lambda^T n\,\widehat{U}(n) e^{2\pi i n\cdot \lambda x},
\end{equation}
where 
\begin{equation}\label{eq:def-directional_gradient}
\nabla_\lambda:=\lambda^T \nabla_Y
\end{equation}
is the directional gradient, $\nabla_Y$ is the full gradient in $\T^N$ and $\widehat{U}(n)$ are the Fourier coefficients of $U$. The main difficulty in estimating the correctors is that the suitable estimates of $\nabla_\lambda U$ (which is easy to obtain for the lifted cell problem) do not give good estimates of $U$ itself, because solving $\widehat{U}(n)$ from $\lambda^T n\,\widehat{U}(n)$ requires dividing by $|\lambda^T n|$ which could be a small divisor.
This motivates us to consider the distribution of $\lambda^T n$ as $|n|$ increases.
For $n\in\Z^N$, we call $n$ the phase frequency and $\lambda^Tn$ its physical frequency. The higher phase frequency of $U(Y)$ does not indicate the higher oscillation of $u(x) = U(\lambda x)$; instead the physical frequency $\lambda^T n$ carries the true frequency information of $u$ in the physical space $\R^d$.

For $\ell\geq 0$, define the index set
\begin{equation*}
\mathcal J_\ell^\lambda:=\left\{j\in\Z:\ 2^{j-1} \le |\lambda^T n| \le 2^{j+1} \text{ for some }n\in\Z^N\text{ with } |n|\le2^{\ell+1}\right\}.
\end{equation*}
The sets $\mathcal J_\ell^\lambda$ record the occupied dyadic rings for the physical frequencies,
which actually contain all the information from the frequency matrix $\lambda$ we need for this paper. A combinatorial counting lemma (see Lemma \ref{lem:band-count}) shows
\begin{equation}\label{eq:band-count}
\#\mathcal J_\ell^\lambda\le C_N (\ell+1)
\end{equation}
uniformly in $\lambda$, where $C_N$ is a constant depending only on $N$.

For $s > 0$, define our key modulus
\begin{equation}\label{eq:def-modulus}
\omega_{\lambda,s}(r):=\sum_{\ell\ge0}2^{-s\ell}\sum_{j\in\mathcal J_\ell^\lambda}\psi(2^jr),\qquad r>0,
\end{equation}
where
\begin{equation}\label{def:psi}
    \psi(t):=\min\{t,t^{-1}\}. 
\end{equation}
The function $\omega_{\lambda,s}$ is continuous, not necessarily monotone, and satisfies the Dini-type condition independent of $\lambda$ (see Lemma \ref{lem:common-modulus})
\begin{equation}\label{eq:dini}
\int_0^\infty\omega_{\lambda,s}(r)\,\frac{\dd r}{r}
=2\sum_{\ell\ge0}2^{-s\ell}\#\mathcal J_\ell^\lambda
\le C_N \sum_{\ell\ge0}(\ell+1)2^{-s\ell}<\infty\,.
\end{equation}
This scale-invariant condition is essential in proving Theorem \ref{thm:main}.
Note that $\omega_{\lambda,s}(r)$ vanishes as $r \to 0$ or $r \to \infty$. However, it may have many peaks when $r = 2^{-j}$ and $j \in \cJ^\lambda_\ell$. The number of these peaks is controlled by \eqref{eq:band-count} and the exponential factor $2^{-s\ell}$ in \eqref{eq:def-modulus} flattens those peaks when $\ell$ is large.

\begin{theorem}[Quantitative sublinear growth]\label{thm:correctors}
Assume \eqref{eq:ellipticity}, \eqref{eq:periodicity} and additionally
\begin{equation}\label{eq:s6-assumptions}
    \|A\|_{C^{m,\alpha}(\T^N)}\le M < \infty,
\end{equation}
where $m\ge0$ is an integer and $0<\alpha\le 1$. Set $s=m+\alpha$. For every $\lambda\in\R^{N\times d}$, the correctors and flux correctors in \eqref{eq:corrector}--\eqref{eq:fluxcorrector} exist, satisfying
\begin{equation}\label{eq:gradientchiphi}
    \| \nabla \chi \|_{L^\infty(\R^d)} + \| \nabla \phi \|_{L^\infty(\R^d)} \le C,
\end{equation}
and, after normalization $\chi(0) = 0$ and $\phi(0) = 0$,
\begin{equation}\label{eq:corrector-growth}
\sup_{|x|\le R}\bigl(|\chi(x)|+|\phi(x)|\bigr)
\le CR\int_R^\infty\omega_{\lambda,s}(r)\,\frac{\dd r}{r}\quad \text{for every }R>0.
\end{equation}
Such normalized correctors and flux correctors that satisfy \eqref{eq:corrector-growth} are unique, and $C$ depends only on $d,N,\Lambda,m,\alpha,M$, and is independent of $\lambda, R$. 
\end{theorem}

Due to \eqref{eq:dini}, the Dini tail in \eqref{eq:corrector-growth} decays to zero and the correctors are sublinear for every fixed $\lambda$. The decay rate can be slow, depending on $\lambda$ via the index sets $\cJ^\lambda_\ell$. Note that the regularity of $A$ can improve the decay rate through the exponent $s$ in \eqref{eq:def-modulus}.

As a corollary, the growth rates of correctors and flux correctors quantify the homogenization scales at which the relative error between a solution and its effective solution becomes small. For example, if $\mathcal{L}_\lambda u = 0$ in $B_{2R}$, then we can find a solution $\bar{u}$ of $-\diver(\overline{A}_\lambda \nabla\bar{u}) = 0$ in $B_R$ such that
\begin{equation*}
    \| u - \bar{u} \|_{\underline L^2(B_R)} \lesssim \int_R^\infty\omega_{\lambda,s}(r)\,\frac{\dd r}{r} \| u  \|_{\underline L^2(B_{2R})},
\end{equation*}
where $\overline{A}_\lambda$ is the homogenized matrix (see \eqref{eq:def-Abar}).
Thus, due to \eqref{eq:dini}, there is a minimal scale such  that the error becomes relatively small, and this scale is a deterministic version of the minimal radius introduced in stochastic homogenization \cite{GNO20}.

\paragraph{Diophantine condition.}
We say $\lambda$ satisfies the Diophantine condition with exponent $\beta > 0$ if for some $c>0$,
\begin{equation}\label{eq:diophantine}
 |\lambda^Tn|\ge c|n|^{-\beta}
 \qquad\text{whenever }n\in\Z^N\text{ and }\lambda^Tn\ne0.
\end{equation}
Using Theorem \ref{thm:correctors}, we show in Theorem \ref{thm:diophantine} that under the Diophantine condition \eqref{eq:diophantine}, for $R\ge2$,
\begin{equation}\label{eq:diophantine-growth0}
\sup_{|x|\le R}\bigl(|\chi(x)|+|\phi(x)|\bigr)\le C
\begin{cases}
1,&\beta<s,\\
\log_2 R,&\beta=s,\\
R^{1-s/\beta},&\beta>s,
\end{cases}
\end{equation}
where $C$ may also depend on $c,\beta$. In particular, $s>\beta$ suffices for bounded correctors and flux correctors. For $0<\alpha<1$ and integer $\beta>s$, we construct two-dimensional scalar equations with diagonal coefficients whose correctors and flux correctors both attain the growth rate $R^{1-s/\beta}$, i.e.,
\begin{equation*}
    \sup_{|x|\le R}|\chi(x)|\simeq \sup_{|x|\le R} |\phi(x)| \simeq R\int_R^\infty\omega_{\lambda,s}(r)\,\frac{\dd r}{r} \simeq R^{1-s/\beta}.
\end{equation*}
This indicates that the sublinear growth rate in \eqref{eq:corrector-growth} is generally optimal. The criterion for bounded correctors in \eqref{eq:diophantine-growth0} under the Diophantine condition is also optimal.

\paragraph{Comparison with earlier work}
We compare Theorem \ref{thm:correctors} (and its consequences) with Shen's \cite{S15} and Armstrong-Gloria-Kuusi's work \cite{AGK16}. They developed a framework for more general almost-periodic coefficients using different moduli to quantify the almost-periodicity. Shen \cite{S15} introduced $\rho(R)$ as in \eqref{eq:Shen-modulus}, and Armstrong-Gloria-Kuusi \cite{AGK16} introduced a higher-order version $\rho_k(R)$ and defined
\begin{equation*}
    \rho_*(R):= \inf_{k\in [1,R]} C^k k! \rho_k(R/k).
\end{equation*}
It was shown that $\rho_*(R) \lesssim R^{-1-\eps}$ implies the boundedness of correctors \cite[Theorem 1.1]{AGK16} and $\rho_*(R) \lesssim R^{-\theta}$ implies their sublinear growth as $R^{1-\theta+\eps}$ for any $\eps>0$ \cite[(1.14)]{AGK16}. In particular, the boundedness of correctors was shown under the Diophantine condition \cite[Section 4]{AGK16}.

In the case of quasiperiodic coefficients, we can compare the modulus $\rho_*$ with our Dini integral. In fact, one can show for $R$ large enough
\begin{equation*}
    \rho_*(R) \le C \int_R^\infty\omega_{\lambda,s}(r)\,\frac{\dd r}{r}.
\end{equation*}
This provides an easy way to estimate $\rho_*(R)$ (not optimal). Recall that $\rho_*(R)$ can decay as $R^{-1-\eps}$. However, the Dini tail decays at most as $R^{-1}$, sufficiently giving the boundedness of correctors. In both \cite{S15} and \cite{AGK16}, under the Diophantine condition, the authors used the Erd\H{o}s-Tur\'{a}n-Koksma inequality to roughly estimate their modulus $\rho$ or $\rho_*$, and hence the smoothness requirements for the coefficients there are not optimal. Our new modulus, depending explicitly on $\lambda$ and $s = m+\alpha$ only, is relatively easy to compute and gives optimal smoothness requirements for coefficients.

\subsection{Proof outline}
We will provide the outline and key insight for the proof of Theorem \ref{thm:main}, while Theorem \ref{thm:correctors} is proved under the same framework. Some unusual techniques will also be briefly explained.

\paragraph{Lipschitz estimate criterion.} We introduce a criterion for the uniform Lipschitz estimate, which could be viewed as a reformulation of \cite[Lemma 8.5]{S17}. The idea is that if we can approximate a solution $u$ at every scale $r>0$ by a solution of an elliptic equation with constant coefficients up to a relative error bounded by $\omega(r)$, and this modulus satisfies the Dini condition $\int_0^\infty \omega(r) r^{-1} \dd r = D_0 < \infty$, then we can conclude that $u$ has the uniform Lipschitz estimate; see Proposition \ref{prop:criterion}. Unlike \cite[Lemma 8.5]{S17}, here we do not require the monotonicity of $\omega(r)$ and even the number of peaks of $\omega(r)$ is not relevant. Only the bound $D_0$ is involved in the final uniform Lipschitz estimate. Eventually, we will show that $\omega(r) = \omega_{\lambda,\kappa}(r)$ as defined in \eqref{eq:def-modulus} for any $\kappa \in (0,\alpha)$. By \eqref{eq:dini}, this will imply the Lipschitz estimate independent of $\lambda$.

\paragraph{Lifted correctors and low-high frequency separation.}

In order to find the constant coefficients of the approximate elliptic operator at every scale $r>0$ and estimate the error, we need to study the correctors and flux correctors. Similarly to \cite{NXZ26}, we begin with the lifted corrector equations with a regularization term
\begin{equation}\label{eq:IntroCorrector}
    \mu^2\X_{\mu,j}-\nlam\cdot\bigl(A(e_j+\nlam\X_{\mu,j})\bigr)=0
    \qquad\text{in }\T^N.
\end{equation}
The positive regularization term guarantees the existence of the solution of \eqref{eq:IntroCorrector} and makes the solutions more regular (depending on $\mu$). Observe that we can use the restriction $\chi_\mu(x) = \X_\mu(\lambda x)$ to recover correctors in physical space $\R^d$. We will call the line $x \mapsto Y_0 + \lambda x$ a physical orbit in $\T^N$ passing through $Y_0 \in \T^N$. Thus, lifting to the torus and restricting to physical orbits are exactly opposite operations in our argument.

To proceed, the general idea in classical homogenization is to write the flux $A(I + \nabla_\lambda \X_\mu)$ as
\begin{equation}\label{eq:low-high-sep}
    A(I + \nabla_\lambda \X_\mu) = \U_\mu + \nabla_\lambda \cdot (\Phi_\mu),
\end{equation}
where $\Phi_\mu$ is the flux corrector tensor, and ideally $\U_\mu$ is expected to be a constant matrix that can be used to define the approximate elliptic operator. However, in our situation with the regularized corrector equation \eqref{eq:IntroCorrector}, $\U_\mu$ need not be constant. The key insight here is that $\U_\mu$ carries mainly low physical frequencies and the high frequencies are enclosed in the flux corrector part $\nabla_\lambda \cdot (\Phi_\mu)$. Therefore, \eqref{eq:low-high-sep} can be viewed as a low-high frequency separation, in the same spirit of scale-separation in multiscale problem; see Subsection \ref{subsec:lowhigh} for a detailed explanation. The issue is that all the frequencies are coupled and the separation is not explicit. Once we accept that $\U_\mu$ mainly contains low physical frequencies, it is natural to further write
\begin{equation*}
    \U_\mu = \U_\mu(0) + \mathcal{R}_\mu,
\end{equation*}
where $\U_\mu(0)$ now is a constant matrix and the remainder $\mathcal{R}_\mu$ is expected to be small along a physical orbit $ \lambda x$ with restricted $|x| < r$. 

\paragraph{Classical and directional Littlewood-Paley decompositions.}

By a formal calculation for flux correctors, the low frequency part $\U_\mu$ in \eqref{eq:low-high-sep} formally satisfies
\begin{equation*}
    (\mu^2 - \Delta_\lambda) \U_\mu = \mu^2 (A(I + \nabla_\lambda \X_\mu) - \nabla_\lambda \X_\mu) =: \mu^2 F,
\end{equation*}
where $\Delta_\lambda$ is the lifting of the Laplace operator in \eqref{eq:def-directional_laplace}; see \eqref{eq:resolvent-cell-flux}. Consequently, we can write
\begin{equation}\label{eq:U-Fourier}
    \widehat{\U}_\mu(n) = \frac{\mu^2}{\mu^2 + 4\pi^2 |\lambda^T n|^2} \widehat{F}(n).
\end{equation}
To show that $\U_\mu$ mainly consists of low physical frequencies, we need to show $\nabla_\lambda \U_\mu$ is small (thus the restriction $\U_\mu(\lambda x)$ varies slowly in a relatively small ball). However, since $\lambda^T n$ appears in the denominator in \eqref{eq:U-Fourier}, we have to balance three parts: size of $\mu$, sizes of $|\lambda^T n|$ and decay of $\widehat{F}(n)$.
This inspires us to introduce two frequency decompositions. The first one is the classical Littlewood-Paley projections $\Delta_\ell$, which localize the phase frequencies with $|n| \simeq 2^\ell$, $\ell \ge 0$. The second one is a directional Littlewood-Paley decomposition with projection $P^\lambda_j$ localizing the physical frequencies with $|\lambda^T n| \simeq 2^j$ for $j \in \Z$. Precisely, there is a smooth function $\rho$ supported in $\{ 1/2 \le |\xi| \le 2\}$ such that for any $V \in L^1(\T^N)$,
\begin{equation*}
    \widehat{P^\lambda_j V}(n) = \rho(2^{-j} \lambda^T n) \widehat{V}(n),
\end{equation*}
and
\begin{equation*}
    V = \sum_{j\in \Z} P^\lambda_j V + \Pi_\lambda V,
\end{equation*}
where $\Pi_\lambda V$ is the projection onto the exact resonant frequencies, i.e., $\widehat{\Pi_\lambda V}(n) = \widehat{V}(n) \one_{ \{ \lambda^T n = 0\} }$. Combined with $\Delta_\ell$ and the definition of the index set $\cJ^\lambda_\ell$, we have
\begin{equation}\label{eq:V-FreqDecomp}
    V = \sum_{\ell \in \N} \sum_{j \in \cJ^\lambda_\ell} P^\lambda_j \Delta_\ell V + \Pi_\lambda V.
\end{equation}
This is because the Fourier multipliers of $\Delta_\ell$ and $P^\lambda_j$ have disjoint supports if $j\notin \cJ^\lambda_\ell$. Note that $\Pi_\lambda V$ remains constant in any physical orbits $Y + \lambda x$. 

We will apply the identity \eqref{eq:V-FreqDecomp} to $\nabla_\lambda \U_\mu$, with $\U_\mu$ given by \eqref{eq:U-Fourier}.
Classically, the LP projection $\Delta_\ell $ carrying the regularity guarantees the decay of $\widehat{F}(n)$. We emphasize that, in order to get the exponential decay of $\Delta_\ell F$, we need the first key ingredient, namely, the uniform $W^{1,p}$ estimate for arbitrary quasiperiodic coefficients in \cite[Corollary 1.2]{NZ24} (see Lemma \ref{lem:physical-cz}).
The second key ingredient, for the directional LP projection, is the counting argument \eqref{eq:band-count} (see Lemma \ref{lem:band-count}), which makes \eqref{eq:V-FreqDecomp} summable. In fact, we can show
\begin{equation*}
    \| P^\lambda_j \Delta_\ell \nabla_\lambda \U_\mu \|_\infty \le C \frac{\mu^2 2^j}{\mu^2 + 2^{2j}} 2^{-\kappa \ell} \le C2^{-\kappa \ell} \mu \psi(2^j \mu^{-1}),
\end{equation*}
where $\psi$ is given by \eqref{def:psi} and $\|\cdot\|_p$ denotes the $L^p$ norm on the torus as in \eqref{eq:pnorm}. At every scale $r>0$, we can take $\mu = r^{-1}$. Summing over $j \in \cJ^\lambda_\ell$ and $\ell \in \N$, together with \eqref{eq:band-count}, gives \eqref{eq:def-modulus} and \eqref{eq:dini}. This argument shows that $\U_\mu$ can have wild oscillation for some scales $r>0$ (precisely, for large $\omega_{\lambda,\kappa}(r)$); however, the bad scales are quantitatively controlled in terms of the Dini condition \eqref{eq:dini}.

We point out that the directional LP decomposition and double frequency decomposition \eqref{eq:V-FreqDecomp} have appeared in literature in the study of quasiperiodic functions and other PDE problems. Particularly, Schippa \cite{Schi25} used exactly the same double decomposition and a lattice counting argument to derive Strichartz estimates for nonlinear Schr\"{o}dinger equations (also see the subsequent work for KdV and Benjamin-Ono equations \cite{Schi26}). The applications of the directional LP projections (possibly in different forms) can also be found in \cite{HTT11,Tsu12,AA24} etc. Nevertheless, we believe this is the first time the above frequency decomposition argument (combined with a combinatorial counting argument) has been used in quasiperiodic homogenization.

\paragraph{Combinatorial counting arguments}

We briefly explain two combinatorial counting arguments applied in this paper. The first one is the counting of index set $\cJ^\lambda_\ell$; see \eqref{eq:band-count} or Corollary \ref{cor:counting}. 

This counting argument has been known and used in different contexts and forms.
Our proof is similar to \cite[p.~477]{Lind65} and \cite[Theorem~2]{LY11}. 
A similar result with a different proof can be found in \cite[Lemma 8.1]{TV09}. 

The second combinatorial counting argument is contained in the proof of Lemma \ref{lem:improve-cancel}, which counts the number of distinct exponents for the power of a trigonometric polynomial. It is used to avoid the exponent loss in the classical Bernstein's inequality (see Lemma \ref{Bern}). We can explain the argument in a simpler form (but in the same spirit). Consider an improved Bernstein's inequality on $\T^N$ (see e.g. \cite[p.~10]{NW78} or \cite[Lemma~3]{Sma98})
\begin{equation*}
    \| V\|_\infty \le |\#\text{spec} (V^m)|^{\frac{1}{2m}} \| V\|_{2m},
\end{equation*}
where $m$ is an integer and $\text{spec}(F)$ is the spectrum of $F$, i.e., the support of $\widehat{F}(n)$. A trivial estimate gives $\#\text{spec} (V^m) \le (\#\text{spec} (V))^m$. The counting argument for a certain class of trigonometric polynomials can result in $\#\text{spec} (V^m) \ll (\#\text{spec} (V))^m$ due to the packing of exponents. This improvement in Lemma \ref{lem:improve-cancel} together with the orbit estimate in Lemma \ref{lem:potentials} allows us to obtain the sharp growth rates in Theorem \ref{thm:correctors} without any loss of exponent $s = m + \alpha$.

\paragraph{Organization of the paper.}
Section~\ref{sec.LipCriterion} proves the Dini criterion for uniform Lipschitz estimate. Section~\ref{sec.LP} introduces the classical and directional LP frequency decompositions and the counting lemma. Section~\ref{sec.corrector} establishes the uniform estimates for (regularized) correctors and flux correctors. Section~\ref{sec.proveLip} proves Theorem \ref{thm:main}. Section~\ref{sec.lineargrowth} proves Theorem \ref{thm:correctors}, establishes explicit growth rates under the Diophantine condition and shows their optimality by examples.

\paragraph{Acknowledgments.} W. N.  is partially supported by NNSF of China (No.12371106) and Scientific Research Project of Anhui Provincial Department of Education (No. 2025AHGXZK20072).  Y. W. and Y. X. are partially supported by NNSF of China (Nos. 12201604, 12371106). J. Z. is partially supported by NNSF of China (Nos. 12494541, 12288201, 12471115).

\paragraph{AI Disclosure Statement.} This project was initiated in May 2026. During the project, the authors had numerous rounds of discussions with various LLMs (DeepSeek-V4-Flash, Kimi K2.6, GPT-5.6 Sol and GPT-6 Astra) to generate original ideas. The first useful progress for a special case in multiscale homogenization setting was obtained with GPT-5.6 Sol in early August. The full uniform Lipschitz estimate in the same setting was obtained with GPT-6 Astra in early September 2026. The authors then guided AI to reformulate the results and their proofs in the current quasiperiodic form with the additional extension on corrector growth rates. The authors carefully verified AI's proofs, carried out substantial simplifications and reconstructions, and then completely rewrote the manuscript. AI also has been used for reference searching and language editing. The authors take full responsibility for the originality and correctness of the paper.

\section{A criterion for Lipschitz estimate}
\label{sec.LipCriterion}

In this section, we prove a criterion for Lipschitz estimate. It states that local approximation by constant-coefficient equations at every scale implies the uniform Lipschitz estimate, provided that the relative error $\omega(r)$ satisfies the Dini-type condition \eqref{eq:dini-functional}. The constant coefficients are allowed to vary at different scales. This criterion has essentially been formulated in \cite[Lemma 8.5]{S17}. We here provide a modified version in a general framework.

Let $\mathcal{P}_1$ denote the set of affine functions on $\R^d$. For a function $u\in L^2(B_r)$, define
\begin{equation*}
    E_1(r; u) = \frac1r \inf_{\ell \in \mathcal{P}_1} \anorm{u - \ell}{B_r}, \qquad E_0(r;u) = \frac1r \inf_{a\in \R} \anorm{u - a}{B_r}.
\end{equation*}

\begin{lemma}\label{lem:Aharmonic}
    Let $A$ be a constant matrix satisfying \eqref{eq:ellipticity}. There exists $C>0$, depending only on $d$ and $\Lambda$, such that if $u$ is $A$-harmonic in $B_r$, then for $\theta \in (0,1/4)$,
    \begin{equation*}
        E_1(\theta r; u) \le C\theta E_1(r; u).
    \end{equation*}
\end{lemma}
\begin{proof}
    The lemma follows from the classical $C^{1,1}$ estimates for second-order elliptic equations with constant coefficients. 
\end{proof}

The main result of this section is the following Lipschitz estimate criterion. 
\begin{proposition}\label{prop:criterion}
Let $A(x)$ be a coefficient matrix on $\R^d$ satisfying the ellipticity and boundedness condition  \eqref{eq:ellipticity}. Suppose that there is a bounded nonnegative continuous function $\omega:(0,\infty) \to \R$, such that
\begin{equation}\label{eq:dini-functional}
\Dini(\omega):=\int_0^\infty\omega(r)\,\frac{\dd r}{r} = D_0 < \infty
\end{equation}
and the following property holds: for any $r>0$ and any $A$-harmonic function $v$ in $B_{4r}$, there exist a uniformly elliptic constant matrix $ \overline{A}_{r}$ and an $\overline{A}_{r}$-harmonic function $\overline{v}_r$ in $B_r$, such that,
\begin{equation}\label{eq:comparison-hypothesis}
\frac1r \anorm{v-\overline{v}_r}{B_r}\le C_*\omega(r)\anorm{\nabla v}{B_{4r}}.
\end{equation}
Let $f\in L^p(B_R)$ and $ g \in C^{0,\gamma}(B_R;\R^d)$ with $p>d$, $0<\gamma\le1$. Then the solution $u\in H^1(B_R)$ to $-\diver(A\nabla u)=f+\diver g$ in $B_R$ satisfies
\begin{equation*}
    \sup_{0<r<R} \anorm{\nabla u}{B_r} \le C\left\{ \anorm{\nabla u}{B_R} + R \pnorm{f}{p}{B_R} + R^\gamma[g]_{C^{0,\gamma}(B_R)} \right\},
\end{equation*}
where the constant $C$ depends on $d, \Lambda, p,\gamma, C_*, D_0$.
\end{proposition}

\begin{proof}
By scaling it is enough to prove the estimate for $R=1$; the rescaling argument is specified at the end of the proof. 

Fix $R_0=1/8$ and set
$$
\alpha_1 = \min\{1-d/p,\gamma\} > 0, \qquad
M = \pnorm{f}{p}{B_1}+[g]_{C^{0,\gamma}(B_1)}.
$$
For each $r\le R_0$, let $w_r \in H^1_0(B_r(x_0))$ be the weak solution of
$$
-\diver( A \nabla w_r)=f+\diver(g - g(0)) \quad \text{in } B_r.
$$

The energy estimate and Poincar\'{e}'s inequality yield
\begin{equation}\label{eq:source-energy}
\begin{split}
\anorm{\nabla w_r}{B_r}+r^{-1}\anorm{w_r}{B_r}
&\le C\{r\pnorm{f}{p}{B_r}+r^\gamma[g]_{C^{0,\gamma}(B_r)}\}\\
&\le Cr^{\alpha_1} M.
\end{split}
\end{equation}
Consider $v_r=u-w_r$, which is $A$-harmonic in $B_r$, i.e., $-\nabla\cdot (A \nabla v_r) = 0$ in $B_r$. The Caccioppoli inequality, after subtracting arbitrary constants, shows
\begin{align*}
    \anorm{\nabla v_r}{B_{r/2}} &\le Cr^{-1} \inf_{a \in \R} \anorm{v_r - a}{B_r} \\
    & \le Cr^{-1} \inf_{a \in \R} \anorm{u - a}{B_r} + Cr^{-1} \anorm{w_r}{B_r} \\
    & \le C\{E_0(r;u)+r^{\alpha_1} M \}.
\end{align*}
By the assumption, there exist a constant matrix $\overline{A}_r$ and an $\overline{A}_r$-harmonic function $\overline{v}_r$ such that
\begin{equation}\label{eq:vr-vrbar}
\begin{aligned}
    \frac1r \anorm{v_r - \overline{v}_r}{B_{r/8}}  & \le C_* \omega(r/8) \anorm{\nabla v_r}{B_{r/2}} \\
    & \le C\omega(r/8) \{E_0(r;u)+r^{\alpha_1} M \}.
\end{aligned}
\end{equation}
Combining \eqref{eq:source-energy} and \eqref{eq:vr-vrbar}, as well as $u = v_r + w_r$, we arrive at
\begin{equation}\label{eq:u-vrbar}
\frac1r\anorm{u-\overline{v}_r}{B_{r/8}}\le C\omega(r/8)\{E_0(r;u)+r^{\alpha_1} M \}+Cr^{\alpha_1} M.
\end{equation}

Now, letting $\theta \in (0,1/32)$ and by the triangle inequality, we have
\begin{equation}\label{eq:E1-u.ttr}
\begin{aligned}
    E_1(\theta r; u) & \le E_1(\theta r; \overline{v}_r) + \frac{1}{\theta r} \anorm{u - \overline{v}_r}{B_{\theta r}} \\
    & \le C\theta E_1(r/8; \overline{v}_r) + \frac{C_\theta }{r} \anorm{u-\overline{v}_r}{B_{r/8}},
\end{aligned}
\end{equation}
where Lemma \ref{lem:Aharmonic} is used in the second line. Moreover, by the triangle inequality again,
\begin{equation}\label{eq:E1-vrbar}
\begin{aligned}
    E_1(r/8; \overline{v}_r) & \le E_1(r/8; u) + \frac{8}{r} \anorm{u-\overline{v}_r}{B_{r/8}} \\
    & \le C E_1(r; u) + \frac{8}{r} \anorm{u-\overline{v}_r}{B_{r/8}}.
\end{aligned}
\end{equation}
Inserting \eqref{eq:E1-vrbar} back into \eqref{eq:E1-u.ttr} and using \eqref{eq:u-vrbar}, we obtain
\begin{equation}\label{eq:E1-u}
    E_1(\theta r; u) \le C\theta E_1(r; u) + C_\theta \omega(r/8) \{E_0(r;u)+r^{\alpha_1} M \} + C_\theta r^{\alpha_1} M.
\end{equation}
Fix $\theta$ small such that $C\theta \le 1/2$ and $\theta^{\alpha_1} \le 1/4$. Define
\begin{equation*}
    H(r) = E_1(r;u) + K_* r^{\alpha_1} M, \qquad \Phi(r) = E_0(r;u) + r^{\alpha_1} M,
\end{equation*}
where $K_*$ is a constant to be determined. Then \eqref{eq:E1-u} is reduced to
\begin{equation*}
    H(\theta r) \le \frac12 E_1(r;u) + C_\theta \omega(r/8) \Phi(r) + (C_\theta + \frac14 K_*) r^{\alpha_1} M.
\end{equation*}
Choosing $K_*$ large such that $C_\theta + \frac14 K_* \le \frac12 K_*$, we obtain
\begin{equation*}
    H(\theta r) \le \frac12 H(r) + C\omega(r/8) \Phi(r).
\end{equation*}
This is \eqref{eq:iter-contract} in Lemma \ref{lem:iteration}.

Next, we verify the rest of the conditions in Lemma \ref{lem:iteration}. Let $\ell_r \in \mathcal{P}_1$ be the affine function achieving the minimum of $E_1(r;u)$ and $h(r) = |\nabla \ell_r| \ge 0$. Then it is easy to see that
\begin{equation*}
    H(r) \le K_* \Phi(r), \qquad \Phi(r) \le H(r) + C h(r), \qquad \sup_{r/2\le t < r} \Phi(t) \le C\Phi(r),
\end{equation*}
and
\begin{equation*}
    |h(\theta r) - h(r)| \le C(H(\theta r) + H(r)).
\end{equation*}
To see the last estimate, we note that
\begin{equation*}
\begin{aligned}
    |h(\theta r) - h(r)| & \le |\nabla \ell_{\theta r} - \nabla \ell_r | \le \frac{C}{\theta r} \anorm{\ell_{\theta r} - \ell_r}{B_{\theta r}} \\
    & \le \frac{C}{\theta r}\anorm{u - \ell_{\theta r}}{B_{\theta r}} + \frac{C}{\theta r}\anorm{u - \ell_{r}}{B_{\theta r}} \\
    & \le C(H(\theta r) + H(r)).
\end{aligned}
\end{equation*}
These constants $C$ are all fixed once $\theta$ is fixed. In addition, $\Dini(\omega(\cdot/8))=\Dini(\omega) = D_0 < \infty$. Hence, all the conditions of Lemma \ref{lem:iteration} are satisfied, which yields that
\begin{equation*}
    \sup_{0<r<1} \Phi(r)  \le C\Phi(1).
\end{equation*}
Furthermore, by the Caccioppoli inequality
\begin{equation*}
    \anorm{\nabla u}{B_r} \le C\Phi(2r)
\end{equation*}
and Poincar\'{e} inequality
\begin{equation*}
    \Phi(r) \le C\anorm{\nabla u}{B_r} + r^{\alpha_1} M,
\end{equation*}
we have
\begin{equation*}
    \sup_{0<r<1/2} \anorm{\nabla u}{B_r} \le C( \anorm{\nabla u}{B_1} + M ).
\end{equation*}
This proves the desired estimate for $R = 1$.

Lastly, we show the scaling argument to generalize the estimate to arbitrary $R>0$. For the original ball $B_R$, set $U(x)=u(Rx)$, which satisfies
$$
-\diver \big(A(Rx) \nabla U \big)=R^2f(Rx)+\diver (R g(Rx)) \qquad \text{in } B_1.
$$
The rescaled equation has comparison modulus $\omega(R\,\cdot)$, and the scale-invariant property of $\Dini(\omega)$ yields $\Dini(\omega(R\,\cdot)) = \Dini(\omega) = D_0$. Moreover,
\begin{align*}
\begin{gathered}
\anorm{\nabla U}{B_1} = R \anorm{\nabla u}{B_R},\\
\pnorm{R^2f(R\,\cdot)}{p}{B_1} =R^2\pnorm{f}{p}{B_R},\\
[Rg(R\,\cdot)]_{C^{0,\gamma}(B_1)} = R^{1+\gamma}[g]_{C^{0,\gamma}(B_R)}.
\end{gathered}
\end{align*}
Consequently, the estimate in $B_R$ for $u$ follows from the corresponding estimate of $U$ in $B_1$. This ends the proof.
\end{proof}

\begin{remark}
    Proposition \ref{prop:criterion} is a full-scale version of Lipschitz estimate criterion. One can also derive a truncated version by replacing the condition \eqref{eq:dini-functional} by
    \begin{equation*}
        \int_\eps^R \omega(r) \frac{\dd r}{r} = D_0 < \infty.
    \end{equation*}
    In this case, we can obtain a large-scale Lipschitz estimate above $\eps$-scale.

    We also point out that Proposition \ref{prop:criterion} applies even to classical periodic homogenization with one scale. In this case, assuming $C^\alpha$ coefficients, we have $\omega(r) = C \min \{(r/\eps)^\alpha, (\eps/r)^\beta \}$, for some $\beta>0$, which satisfies \eqref{eq:dini-functional}. Note that $(\eps/r)^\beta$ is due to the convergence rate for $r > \eps$ and $(r/\eps)^\alpha$ is due to the regularity of coefficients for $r< \eps$. It is easy to verify that such $\omega(r)$ satisfies the Dini condition \eqref{eq:dini-functional} with only one peak at $r = \eps$. Similarly for multiscale problems under \eqref{eq:multiscale-A}, one will have at most $n$ peaks at $r = \eps_i$.
\end{remark}

The following is a version of iteration lemma similar to \cite[Lemma 8.5]{S17}. The difference in our context is that we do not require the modulus $\omega(r)$ to be monotonic and the Dini-type integral \eqref{eq:dini-functional} is defined in the entire $(0,\infty)$ (hence across all scales).
\begin{lemma}\label{lem:iteration}
Let $\omega:(0,\infty) \to [0,\infty)$ be a bounded nonnegative continuous function as in \eqref{eq:dini-functional} and $R_0 > 0$. Let $H,\Phi,h$ be nonnegative continuous functions on $(0,R_0]$. Suppose that, for some $0< \theta, \delta <1$ and $C_0 \ge 1$, and any $r\in (0,R_0]$,
\begin{align}
H(\theta r)&\le \delta H(r)+C_0\omega(r)\Phi(r),\label{eq:iter-contract}\\
|h(\theta r)-h(r)|&\le C_0\{H(r)+H(\theta r)\},\label{eq:iter-slope}\\
C_0^{-1} \{ H(r)+h(r) \} &\le \Phi(r)\le C_0\{H(r)+h(r)\},\label{eq:iter-equivalence} \\
\sup_{r/2 \le t\le r} \Phi(t) & \le C_0 \Phi(r). \label{eq:doubling}
\end{align}
If $\Dini(\omega)\le D_0<\infty$, then
\begin{equation}\label{eq:iter-conclusion}
\sup_{0<r\le R_0}\Phi(r)+\int_0^{R_0}H(r)\,\frac{\dd r}{r}
\le C\Phi(R_0),
\end{equation}
where $C$ depends only on $\delta, C_0, \theta,D_0$.
\end{lemma}
\begin{proof}
We first claim that there exists $\rho \in [\theta R_0, R_0]$ such that
\begin{equation}\label{eq:DiscreteDini}
    \sum_{j=0}^\infty \omega(\theta^j \rho) \le C_\theta \int_0^{R_0} \omega(r) \frac{\dd r}{r}.
\end{equation}
To see this, note that
\begin{equation*}
    \frac{1}{(1-\theta) R_0} \int_{\theta R_0}^{R_0} \sum_{j=0}^\infty \omega(\theta^j \rho) \dd \rho \le C_\theta \int_{\theta R_0}^{R_0} \sum_{j=0}^\infty \omega(\theta^j \rho) \frac{\dd \rho}{\rho} = C_\theta \int_0^{R_0}  \omega(r) \frac{\dd r}{r},
\end{equation*}
which, by the mean value theorem, implies the existence of $\rho \in [\theta R_0, R_0]$ satisfying \eqref{eq:DiscreteDini}.

Write $H_j=H(\theta^j\rho), \Phi_j = \Phi(\theta^j \rho)$, $h_j=h(\theta^j\rho)$, and $\omega_j=\omega(\theta^j\rho)$. To compensate for the possible increase of $h_j$ through the contraction of $H_j$, we combine $h_j$ with a multiple of $H_j$. In detail, define
$$S_j=h_j+K_0 H_j, $$
where $K_0 > 0$ is a constant to be determined.
Then applying \eqref{eq:iter-contract}, \eqref{eq:iter-slope} and \eqref{eq:iter-equivalence} in suitable order gives
\begin{equation*}
    \begin{aligned}
        S_{j+1} & = h_{j+1} + K_0 H_{j+1} \\
        & \le h_j + |h_{j+1} - h_j| + K_0 \delta H_j + K_0 C_0 \omega_j \Phi_j \\
        & \le h_j + C_0 (H_{j+1} +  H_j) + K_0 \delta H_j + K_0 C_0 \omega_j \Phi_j \\
        & \le h_j + (C_0 + K_0 \delta) H_j + C_0(\delta H_j + C_0 \omega_j \Phi_j) + K_0 C_0 \omega_j \Phi_j \\
        & \le h_j + (C_0 + K_0 \delta + C_0\delta )H_j + (C_0^2 + C_0 K_0) \omega_j C_0 (H_j + h_j) \\
        & \le (1 + (C_0^3 + C_0^2 K_0) \omega_j ) h_j + (C_0 + K_0 \delta + C_0 \delta + (C_0^3 + C_0^2 K_0) \omega_j ) H_j.
    \end{aligned}
\end{equation*}
Choose $K_0$ large enough such that $C_0 + K_0\delta + C_0 \delta \le K_0$, which is possible since $\delta \in (0,1)$. Set $C_1 = C_0^3 + C_0^2 K_0$. Then it follows that
\begin{equation*}
    S_{j+1} \le (1 + C_1 \omega_j ) h_j + (1 + C_1 \omega_j) K_0 H_j = (1 + C_1 \omega_j) S_j.
\end{equation*}
Consequently, for every integer $J\ge1$,
$$
S_J\le S_0\prod_{j=0}^{J-1}(1+C_1 \omega_j)
\le S_0\exp \bigg(C_1\sum_{j=0}^{J-1}\omega_j\bigg)
\le S_0\exp( C_\theta D_0).
$$
Here we used $1+t\le e^t$ for $t\ge0$ and \eqref{eq:DiscreteDini}.

By \eqref{eq:iter-equivalence}, $S_0 \le C_0 K_0 \Phi(\rho)$, and hence
\begin{equation}\label{eq:Discrite-Phi}
    \sup_{j \ge 1} \Phi(\theta^j \rho) \le C_0 \sup_{j\ge 1} S_j \le C_0^2 K_0 \exp( C_\theta D_0) \Phi(\rho). 
\end{equation}
By applying \eqref{eq:doubling} finite times, for any $r\in (0,R_0]$, we have
\begin{equation}\label{eq:thetaing}
    \sup_{\theta r \le t < r} \Phi(t) \le C_\theta \Phi(r). 
\end{equation}
This, combined with \eqref{eq:Discrite-Phi},  proves the first part of \eqref{eq:iter-conclusion}.

To see the second part, integrating \eqref{eq:iter-contract} against $r^{-1}$ over $(0,R_0)$, we have
\begin{equation*}
    \int_0^{R_0} H(\theta r) \frac{\dd r}{r} \le \delta \int_0^{R_0} H(r) \frac{\dd r}{r} + C_0 \int_0^{R_0} \omega(r) \frac{\dd r}{r} \sup_{0<r < R_0} \Phi(r).
\end{equation*}
Hence, by \eqref{eq:dini-functional} and the first part of \eqref{eq:iter-conclusion},
\begin{equation*}
    (1-\delta) \int_0^{\theta R_0} H(r) \frac{\dd r}{r} \le \delta \int_{\theta R_0}^{R_0} H(r) \frac{\dd r}{r} + C \Phi(R_0).
\end{equation*}
This, combined with \eqref{eq:iter-equivalence} and \eqref{eq:thetaing}, gives
\begin{equation*}
    \int_0^{R_0} H(r) \frac{\dd r}{r} \le C \Phi(R_0).
\end{equation*}
This proves the second part of \eqref{eq:iter-conclusion}.
\end{proof}

\section{Frequency decompositions and combinatorial counting argument}
\label{sec.LP}

In this section, we introduce two frequency decompositions, the Littlewood--Paley decomposition in the phase frequency and the directional Littlewood-Paley decomposition in the physical frequency. 

For $V\in L^1(\T^N)$, recall that its Fourier coefficients are given by
\begin{equation}\label{eq:phase-fourier}
\widehat V(n)=\int_{\T^N}V(Y)e^{-2\pi i n\cdot Y}\dd Y,
\qquad n\in\Z^N.
\end{equation}

Let $\vartheta$ be a fixed smooth cutoff function on $\R$ such that $\vartheta(s)=1$ for $|s|\le1$ and $\vartheta(s)=0$ for $|s|\ge2$. We define the nonhomogeneous Littlewood-Paley operators/projections to periodic functions on $\T^N$ by
\begin{equation}\label{eq:phase-blocks}
\widehat{\Delta_0V}(n):=\vartheta(|n|)\widehat V(n),\qquad
\widehat{\Delta_\ell V}(n)
:=(\vartheta(2^{-\ell}|n|)-\vartheta(2^{1-\ell}|n|))\widehat V(n),~~\ell\ge1.
\end{equation}
It follows that
\begin{equation}\label{eq:phase-support}
\begin{aligned}
\supp\widehat{\Delta_0V}
&\subseteq\{n\in\Z^N:|n|\le2\},\\
\supp\widehat{\Delta_\ell V}
&\subseteq\{n\in\Z^N:2^{\ell-1}\le|n|\le2^{\ell+1}\},
\qquad \ell\ge1.
\end{aligned}
\end{equation}

Next, set
\begin{equation*}
\rho(\xi):=\vartheta(|\xi|)-\vartheta(2|\xi|)\qquad \text{for }\xi\in \R^d,
\end{equation*}
which satisfies $\supp\rho\subseteq\{1/2\le|\xi|\le2\}$. We define the homogeneous directional Littlewood-Paley projections by
\begin{equation}\label{eq:physical-bands}
\widehat{P_j^\lambda V}(n)
:=\rho(2^{-j}\xi_\lambda(n))\widehat V(n),
\qquad n\in\Z^N,\quad j\in\Z,
\end{equation}
where 
\begin{equation}\label{eq:physical-frequency}
\xi_\lambda(n):=\lambda^T n\in\R^d
\end{equation}
is the physical frequency in the direction $\lambda$ associated with the phase frequency $n$.

Besides, the resonance projection $\Pi_\lambda$ is defined by
\begin{equation*}
    \widehat{\Pi_\lambda V}(n):= \widehat{V}(n) \one_{\{\xi_\lambda(n) = 0 \}},
\end{equation*}
or equivalently, in the sense of distributions
\begin{align}\label{def_Pi_lambda}
    \Pi_\lambda V(Y) = \sum_{\xi_\lambda(n) = 0} \widehat{V}(n) e^{2\pi i n\cdot Y}.
\end{align}

These operators are all Fourier multiplier operators and commute with each other. Note that the Littlewood--Paley operators $\Delta_\ell$ are independent of $\lambda$, while the directional projections depend on $\lambda$ through $\xi_\lambda$. 

Now, for any $V\in L^1(\T^N)$, we have, in the sense of tempered distributions,
\begin{equation}\label{eq:full-decomposition}
     V = \sum_{j\in \Z} P^\lambda_{j} V + \Pi_\lambda V=\sum_{\ell\in \N}\sum_{j\in \Z} P^\lambda_{j}\Delta_\ell V + \Pi_\lambda V.
\end{equation}
We point out that the double summation in \eqref{eq:full-decomposition} contains many empty terms. In fact, for $\ell\ge0$, define the index set
\begin{equation}\label{eq:occupied-bands}
\cJ_\ell^\lambda:=
\big\{j\in\Z: 1/2\leq 2^{-j}|\xi_\lambda(n)|\leq 2
\text{ for some }n\in\Z^N,\ |n|\le2^{\ell+1}\big\}.
\end{equation}
By \eqref{eq:phase-support} and the definitions of
$\cJ_\ell^\lambda$ and $P_j^\lambda$,
$P_j^\lambda\Delta_\ell V=0$ for $j\notin\cJ_\ell^\lambda$. 
Therefore, \eqref{eq:full-decomposition} can be rewritten as
\begin{equation}\label{eq:lp-reconstruction}
 V = \sum_{\ell \in \N} \sum_{j\in\cJ_\ell^\lambda}P_j^\lambda\Delta_\ell V + \Pi_\lambda V.
\end{equation}
The following lemma estimates the number of occupied nonresonance-frequency bands. 

\begin{lemma}\label{lem:band-count}
Let $T:\R^N\to\R^d$ be a linear map and $R\ge2$.
Then the set
\begin{equation*}
\{|Tn|:n\in\Z^N,\ |n|_\infty\le R,\ Tn\ne0\}
\end{equation*}
is contained in the union of at most $C_N\log R$ dyadic intervals $[2^j,2^{j+1})$,
where $j\in\Z$, $|n|_\infty=\max_{1\le i\le N}|n_i|$ and the constant $C_N$ depends only on $N$.
\end{lemma}
\begin{proof}

Let
\begin{equation*}
\mathcal I
=\big\{j\in \Z:\text{ there exists } n\in\Z^N\text{ such that } |n|_\infty\le R,\
2^j\le|Tn|<2^{j+1}\big\}.
\end{equation*}
To show the lemma, we only need to count the number of elements in $\mathcal{I}$.

Let $\mathcal I_{\rm even} = \mathcal I \cap 2\Z$ and $\mathcal I_{\rm odd} = \mathcal I \cap (2\Z + 1)$. The estimates of $\# \mathcal I_{\rm even}$ and $\# \mathcal I_{\rm odd}$ are similar and it suffices to estimate $\# \mathcal I_{\rm even}$.

Assume $\#\mathcal I_{\rm even} =m \ge1$, and write
\begin{equation*}
\mathcal I_{\rm even}=\{j_1,\ldots,j_m\},
\qquad j_1<\cdots<j_m.
\end{equation*}
By definition, there are $n_1,\ldots,n_m\in\Z^N$ such that
\begin{equation*}
|n_i|_\infty\le R,\qquad
2^{j_i}\le|Tn_i|<2^{j_i+1},
\qquad 1\le i\le m.
\end{equation*}
Since $j_i$'s are even numbers, we have
\begin{equation*}
j_{i+1}\ge j_i+2,\qquad
2|Tn_i|<2^{j_i+2}\le2^{j_{i+1}}\le|Tn_{i+1}|
\end{equation*}
for $1\le i<m$. It follows that, for any $1\le i_0\le m$,
\begin{equation}\label{eq:separated-representatives}
\sum_{i=1}^{i_0-1}|Tn_i|
\le |Tn_{i_0}|\sum_{i=1}^{i_0-1}2^{i-i_0}
=(1-2^{1-i_0})|Tn_{i_0}|<|Tn_{i_0}|.
\end{equation}

Now we claim the elements in 
\begin{equation*}
    \bigg\{ T\Big(\sum_{i=1}^m \eps_i n_i \Big): \eps = (\eps_i)_{i=1}^m \in \{0,1 \}^m \bigg\}
\end{equation*}
are pairwise distinct. To see this, let $\varepsilon,\varepsilon'\in\{0,1\}^m$ with
$\varepsilon\ne\varepsilon'$, and let 
\begin{equation*}
i_0 = \max\{1\le i\le m:\varepsilon_i\neq \varepsilon_i'\}.
\end{equation*}
By \eqref{eq:separated-representatives},
\begin{equation*}
\left|T\left(\sum_{i=1}^m
(\varepsilon_i-\varepsilon_i')n_i\right)\right|
\ge |\varepsilon_{i_0}-\varepsilon_{i_0}'|\,|Tn_{i_0}|
-\sum_{i=1}^{i_0-1}|\varepsilon_i-\varepsilon_i'|\,|Tn_i|
\ge |Tn_{i_0}|-\sum_{i=1}^{i_0-1}|Tn_i|>0.
\end{equation*}
This proves the claim.

Since $T$ is a linear map, the above claim implies that the $2^m$ vectors $\{ \sum_{i=1}^m\varepsilon_i n_i \}$ are pairwise distinct for $\eps \in \{ 0,1\}^m$. Moreover,
\begin{equation*}
\sum_{i=1}^m\varepsilon_i n_i\in\Z^N,\qquad
\left|\sum_{i=1}^m\varepsilon_i n_i\right|_\infty
\le\sum_{i=1}^m|n_i|_\infty\le mR.
\end{equation*}
Consequently,
\begin{equation*}
  2^m=\#\left\{\sum_{i=1}^m\varepsilon_i n_i:
    \varepsilon\in\{0,1\}^m\right\}\le\#\big(\Z^N\cap[-mR,mR]^N\big)\le(3mR)^N,
\end{equation*}
which yields
\begin{equation*}
\begin{aligned}
m
&\le N\log_2(3R)+N\log_2m\\
&\le N\log_2(3R)
  +N\left(\log_2(2N)+\frac{m}{2N}\right)\\
&=N\log_2(6NR)+\frac{m}{2}.
\end{aligned}
\end{equation*}
Hence, 
\begin{equation*}
\#\mathcal I_{\rm even}=m\le2N\log_2(6NR).
\end{equation*}
Combined with the same estimate for $\mathcal I_{\rm odd}$, we arrive at
\begin{equation*}
\#(\mathcal I_{\rm even}\cup\mathcal I_{\rm odd})
\le4N\log_2(6NR)
\le4N\big(1+\log_2(6N)\big)\log_2R
\le C_N\log R,
\end{equation*}
as $\log_2R\ge1$ for $R\ge2$. The proof is complete. 
\end{proof}

\begin{corollary}\label{cor:counting}
  For each $\lambda \in \mathbb{R}^{N\times d}$ and $\ell\geq 0$, we have
  \begin{equation}\label{eq:band-cardinality}
\#\cJ_\ell^\lambda\le C_N(\ell+1),
\end{equation}
where $C_N$ depends only on $N$.
\end{corollary}
\begin{proof}
  Let
\begin{equation*}
\mathcal I_\ell^\lambda=\big\{s\in\Z: 2^s\le|\xi_\lambda(n)|<2^{s+1}\text{ for some }n\in \Z^N,\ |n|\le2^{\ell+1}\big\}.
\end{equation*}
For $j\in\cJ_\ell^\lambda$, there exists $n\in\Z^N$ such that
\begin{equation*}
|n|\le2^{\ell+1},\qquad
2^{j-1}\le|\xi_\lambda(n)|\le2^{j+1}.
\end{equation*}
Then one can find the unique $s\in \{j-1, j, j+1\}$ such that
\begin{equation*}
2^s\le|\xi_\lambda(n)|<2^{s+1}.
\end{equation*}
This implies that $s\in\mathcal I_\ell^\lambda$. Comparing the definitions of $\cJ^\lambda_\ell$ and $\mathcal{I}^\lambda_\ell$, we know 
\begin{equation*}
\cJ_\ell^\lambda\subseteq\bigcup_{s\in\mathcal I_\ell^\lambda}\{s-1,s,s+1\}.
\end{equation*}
Since $|n|_\infty\le|n|$, we may apply Lemma \ref{lem:band-count} with
$Tn=\lambda^Tn$ and $R=2^{\ell+1}$ to obtain
\begin{equation*}
\#\cJ_\ell^\lambda
\le3\,\#\mathcal I_\ell^\lambda
\le3C_N\log(2^{\ell+1})
\le C_N(\ell+1),
\end{equation*}
where $C_N$ depends only on $N$.
\end{proof}

\begin{remark}
Without Lemma \ref{lem:band-count}, one may only obtain a rough bound for $\cJ_\ell^\lambda$. For each $n\in\Z^N$, $1/2\leq 2^{-j}|\xi_\lambda(n)|\leq 2$ holds for at most three indices $j$. It follows that
\begin{equation*}
\#\cJ_\ell^\lambda
\le3\,\#\{n\in\Z^N:|n|\le2^{\ell+1}\}
\le C_N2^{N\ell},
\end{equation*}
whereas \eqref{eq:band-cardinality} gives the linear bound $\#\cJ_\ell^\lambda\le C_N(\ell+1)$.
\end{remark}

We end this section with some useful estimates for Littlewood--Paley projections and identifying the resonance projection with an orbit average. Throughout the paper, for $\Z^N$-periodic function $V$, we write
\begin{equation}\label{eq:pnorm}
    \| V \|_p := \| V \|_{L^p(\T^N)}
\end{equation}
for any $1\le p\le \infty$.

\begin{lemma}[Bernstein's inequality \cite{BCD11}]\label{Bern}
For $V\in L^1(\T^N)$, $k, \ell\in\N$ and $1\le q<\infty$,
\begin{equation*}
\|\nabla^k\Delta_\ell V\|_{\infty}
\le C^{k+1}2^{(k+N/q)\ell}\|\Delta_\ell V\|_{q}\,,
\end{equation*}
where $C$ depends only on $N$.
\end{lemma}

\begin{lemma}\label{diff-block}
  Let $V\in L^q(\mathbb{T}^N)$ for some $1\leq q\leq \infty$. If $0<\alpha\le1$ and 
\begin{equation*}
\|D_ZV\|_{q}\le C_0|Z|^\alpha\quad \text{for all } Z\in\R^N
\end{equation*}
with $D_ZV(Y)=V(Z+Y)-V(Y)$, then for $\ell\ge1$
\begin{equation*}
\|\Delta_\ell V\|_{q}\le C2^{-\alpha\ell},
\end{equation*}
where $C$ depends only on $N, \alpha, C_0$ and the cutoff $\vartheta$. 
\end{lemma}

\begin{proof}
The proof is also standard and we provide one for the reader's convenience.
For $\ell \geq 1$, denote
\begin{equation}\label{eq:kernel}
\begin{gathered}
\mathcal K(Z):=\int_{\R^N}e^{2\pi iZ\cdot\zeta}
[\vartheta(|\zeta|)-\vartheta(2|\zeta|)]\dd \zeta,\quad \text{and}\quad
\mathcal K_\ell(Z):=2^{N\ell}\mathcal K(2^\ell Z).
\end{gathered}
\end{equation}
Evidently, by the Fourier transform,
\begin{equation*}
\widehat{\mathcal K_\ell}(\zeta)
=\vartheta(2^{-\ell}|\zeta|)-\vartheta(2^{1-\ell}|\zeta|),
\end{equation*}
and, for each multi-index $\beta\in\N^N$ and $\ell\geq 1$, 
\begin{equation}\label{eq:phase-kernel-moments}
\int_{\R^N}Z^\beta\mathcal K_\ell(Z)\dd Z=(-2\pi i)^{-|\beta|}(\partial^\beta\widehat{\mathcal K_\ell})(0)=0.
\end{equation}
The operator $\Delta_\ell$ acts on the periodic extension of $V$ by convolution with $\mathcal{K}_\ell$. Therefore,
\begin{equation*}
\Delta_\ell V(Y)
=\int_{\R^N}\mathcal K_\ell(Z)V(Y-Z)\dd Z
=\int_{\R^N}\mathcal K_\ell(Z)[V(Y-Z)-V(Y)]\dd Z.
\end{equation*}
By Minkowski's inequality and the assumption of the lemma, this implies
\begin{equation*}
\|\Delta_\ell V\|_{q}
\le\int_{\R^N}|\mathcal K_\ell(Z)|
\|D_{-Z}V\|_{q}\dd Z
\le C_0\int_{\R^N}|\mathcal K_\ell(Z)||Z|^\alpha\dd Z
\le C2^{-\alpha\ell},
\end{equation*}
where in the last step we have used the fact
\begin{equation}\label{eq:phase-kernel-cancellation}
\int_{\R^N}|\mathcal K_\ell(Z)||Z|^t\dd Z
=2^{-t\ell}\int_{\R^N}|\mathcal K(W)||W|^t\dd W
\le C_t2^{-t\ell}\quad \text{for any }t\geq 0.
\end{equation}
The proof is complete.
\end{proof}

\begin{lemma}\label{lem:orbit-averages}
For $\lambda\in\R^{N\times d}$ and $V\in C(\T^N)$, the following statements hold.
\begin{enumerate}[(i)]
\item The projection $\Pi_\lambda V$ has a continuous representative given by the orbit average
\begin{equation}\label{eq:resonant-orbit-average}
(\Pi_\lambda V)(Y)=\lim_{R\to\infty}\fint_{B_R}V(Y+\lambda x)\dd x.
\end{equation}
The convergence is uniform in $Y$ for each fixed $\lambda$.
\item (Monotonicity) If $V_1\geq V_2$ on $\T^N$, then $\Pi_\lambda V_1\geq \Pi_\lambda V_2$ on $\T^N$.
\item (Orbit invariance) For all $Y\in\T^N$ and $x\in\R^d$,
\begin{equation}\label{eq:resonant-orbit-invariance}
(\Pi_\lambda V)(Y+\lambda x)=(\Pi_\lambda V)(Y).
\end{equation}
\item $\|V\|_{1}\le\sup_{Y\in\T^N}(\Pi_\lambda|V|)(Y)$.
\end{enumerate}
\end{lemma}

\begin{proof}
  For $\xi_\lambda(n)\ne0$, it follows from the Riemann-Lebesgue Lemma that
\begin{equation*}
\begin{aligned}
\fint_{B_R}e^{2\pi i\xi_\lambda(n)\cdot x}\dd x
&=|B_1|^{-1}\int_{B_1} e^{2\pi i R\xi_\lambda(n)\cdot x}\dd x\longrightarrow0,\quad \text{as }R\rightarrow \infty. 
\end{aligned}
\end{equation*}
For $\xi_\lambda(n)=0$, the average is $1$. Noticing that
\begin{align*}
  e^{2\pi i n\cdot[Y+\lambda x]}=e^{2\pi i n\cdot Y}\cdot e^{2\pi i\xi_\lambda(n)\cdot x},
\end{align*}
we immediately know \eqref{eq:resonant-orbit-average} holds for any trigonometric polynomials. Since the set of trigonometric polynomials is a dense subspace of $C(\T^N)$ and the linear operator $\Pi_\lambda$ on this subspace has norm at most $1$, one can see that \eqref{eq:resonant-orbit-average} holds for $V\in C(\mathbb{T}^N)$, uniformly in $Y$.

The item (ii) is a direct consequence of (i), and (iii) follows from \eqref{def_Pi_lambda}.

Lastly, since by \eqref{def_Pi_lambda}
\begin{align*}
  \int_{\T^N} V(Y)\dd Y = \widehat{V}(0) =\int_{\T^N} \Pi_\lambda V(Y)\dd Y,
\end{align*}
we have
\begin{align}\label{es:L1Linfty}
\|V\|_1=\int_{\T^N}(\Pi_\lambda|V|)(Y)\dd Y
\le\sup_{Y\in\T^N}(\Pi_\lambda|V|)(Y),
\end{align}
which proves (iv).
\end{proof}

\section{Uniform estimates of correctors and flux correctors}
\label{sec.corrector}
In this section we establish various estimates for correctors and flux correctors, and scale-adapted low-high frequency separation.

\subsection{Correctors}

As in the previous study on quasiperiodic homogenization, the true corrector of operator \eqref{eq:operator} is given by
\begin{equation*}
    -\diver (A_\lambda\nabla \chi_j)=\diver(A_\lambda e_j)\qquad \text{in }\R^d. 
\end{equation*}
Due to the quasiperiodic structure of the coefficient, $\chi_j$ would has also the form
\begin{align*}
    \chi_j(x)=\X_j(\lambda x),
\end{align*}
and $\X_j = \X_j(Y):\T^N \to \R$ satisfies the lifted cell problem
\begin{align*}
    -\nabla_\lambda\cdot(A\nabla_\lambda\X_j)=\nabla_\lambda\cdot(Ae_j)\qquad \text{in }\T^N.
\end{align*}
where $\nabla_\lambda$ is defined by \eqref{eq:def-directional_gradient}. However, this lifted equation may be degenerate and the control of $\X_j$ itself is absent. Therefore, as in \cite{NXZ26}, we recover partial coercivity by adding a positive zeroth-order term. Precisely, we consider instead the equation of approximate correctors
\begin{equation}\label{eq:cell-equation}
\mu^2\X_{\mu,j}-\nlam\cdot\bigl(A(e_j+\nlam\X_{\mu,j})\bigr)=0
\qquad\text{in }\T^N,
\end{equation}
where $\mu>0$ is a parameter. 

Before studying the correctors, we present a Calder\'on--Zygmund estimate for equations with quasiperiodic coefficients from \cite{NZ24}.

\begin{lemma}\label{lem:physical-cz}
Suppose that $A$ satisfies \eqref{eq:ellipticity}--\eqref{eq:holder} and $\lambda\in \mathbb{R}^{N\times d}$. Let $R>0$ and $H\in L^q(B_{5R};\R^d)$ with $2<q<\infty$. Suppose that $w\in H^1(B_{5R})$ is a weak solution of
\begin{equation*}
-\diver(A_\lambda\nabla w)=\diver H\quad\text{in }B_{5R}.
\end{equation*}
Then
\begin{equation*}
\pnorm{\nabla w}{q}{B_R}
\le C_q\big\{\anorm{\nabla w}{B_{5R}}
+\pnorm{H}{q}{B_{5R}}\big\},
\end{equation*}
where $C_q$ depends only on $q$ and the character of $A$.
\end{lemma}

\begin{proof}
  This estimate is a quasiperiodic version of  \cite[Theorem 4.1]{NZ24} or a localized version of \cite[Corollary 1.2]{NZ24}. It can be derived following the proof of \cite[Corollary 1.2]{NZ24}. For completeness, we present the argument here.

  For $1\leq i\leq N$, $1\leq j\leq d$, set
  \begin{align*}
    \varepsilon_{ij}=
    \begin{cases}
      |\lambda_{ij}|^{-1} &\lambda_{ij}\neq 0,\\
      1 & \lambda_{ij}= 0. 
    \end{cases}
  \end{align*}
  Write
  \begin{align*}
    A(\lambda x)&=A\bigg(\sum_{j=1}^d\lambda_{1j}e_j\cdot x, \cdots, \sum_{j=1}^d\lambda_{Nj}e_j\cdot x\bigg)\\&=A\bigg(\sum_{j=1}^d\mathrm{sgn}(\lambda_{1j})e_j\cdot x/\varepsilon_{1j}, \cdots, \sum_{j=1}^d\mathrm{sgn}(\lambda_{Nj})e_j\cdot x/\varepsilon_{Nj}\bigg),
  \end{align*}
  where $e_j$ is the $j$th coordinate vector in $\R^d$ and $$\mathrm{sgn}(s)=\begin{cases}
      s/|s| & \text{if }s\neq 0,\\
      0 & \text{if }s=0.
  \end{cases}$$ By setting
  \begin{align*}
    B(y_{11}, \cdots, y_{ij}, \cdots, y_{Nd})=A\bigg(\sum_{j=1}^d\mathrm{sgn}(\lambda_{1j})e_j\cdot y_{1j}, \cdots, \sum_{j=1}^d\mathrm{sgn}(\lambda_{Nj})e_j\cdot y_{Nj}\bigg)
  \end{align*}
  for independent variables $y_{ij}\in\R^d$, we have
  \begin{align*}
    A(\lambda x)=B(x/\varepsilon_{11}, \cdots, x/\varepsilon_{ij}, \cdots, x/\varepsilon_{Nd}). 
  \end{align*}
  It is easy to see that $B$ is $\Z^d$-periodic in each $y_{ij}\in\mathbb{R}^d$, and satisfies the same ellipticity and H\"{o}lder conditions as $A$. Thus, $w$ is a solution to
  \begin{align*}
    -\diver(B(x/\varepsilon_{11}, \cdots, x/\varepsilon_{ij}, \cdots, x/\varepsilon_{Nd})\nabla w)=\diver H\quad \text{in }B_{5R}, 
  \end{align*}
  which is a multiscale equation of $Nd$ scales. By \cite[Theorem 4.1]{NZ24}, we obtain the desired estimate immediately.
\end{proof}

To deal with equation \eqref{eq:cell-equation}, set
\begin{align*}
  H_\lambda^1(\T^N)=\{v\in L^2(\T^N):\nlam v\in L^2(\T^N;\R^d)\}
\end{align*}
endowed with the norm $\|v\|_{H^1_\lambda(\mathbb{T}^N)}=(\|v\|_2^2+\|\nabla_\lambda v\|_2^2)^{1/2}$, where derivatives are taken in the sense of distributions. It is essential to verify that $H^1_\lambda(\mathbb{T}^N)$ is a real Hilbert space. For $\mu>0$, it has the equivalent norm
\begin{align}\label{norm_H1lambda_equiv}
  (\mu^2\|v\|_2^2+\|\nlam v\|_2^2)^{1/2}. 
\end{align}

\begin{lemma}\label{lem:torus-estimate}
  Let $\mu>0$, $2\leq q<\infty$, and $H\in L^q(\T^N;\R^d)$. Then there exists a unique weak solution $w\in H^1_\lambda(\T^N)$ to
\begin{equation}\label{eq:torus-equation}
\mu^{2}w-\nlam\cdot(A\nlam w)=\nlam\cdot H\quad\text{in }\T^N,
\end{equation}
satisfying  
\begin{equation}\label{eq:torus-estimate}
\mu\|w\|_{q}+\|\nlam w\|_{q}\le C_q\|H\|_{q},
\end{equation}
where $C_q$ depends only on $q$ and the character of $A$. 
\end{lemma}
\begin{proof}
  The existence and uniqueness of the weak solution $w\in H^1_\lambda(\T^N)$ follow from the Lax-Milgram theorem. Indeed, we consider the bilinear form
  \begin{align*}
    a_\mu(v,w)=\mu^2\int_{\T^N}vw+\int_{\T^N}A\nlam v\cdot\nlam w
  \end{align*}
on $H^1_\lambda(\mathbb{T}^N)$. By using the equivalent norm \eqref{norm_H1lambda_equiv} denoted still by $\|\cdot\|_{H^1_\lambda(\mathbb{T}^N)}$, one can see
\begin{align*}
  a_\mu(v, w)\leq \Lambda\|v\|_{H^1_\lambda(\mathbb{T}^N)}\|w\|_{H^1_\lambda(\mathbb{T}^N)},\qquad
  a_\mu(v, v)\geq \Lambda^{-1}\|v\|_{H^1_\lambda(\mathbb{T}^N)}^2,
\end{align*}
i.e., $a_\mu$ is bounded and coercive on $H^1_\lambda(\mathbb{T}^N)$. By the Lax-Milgram theorem, equation \eqref{eq:torus-equation} admits a unique weak solution $w$ in $H^1_\lambda(\mathbb{T}^N)$ and \eqref{eq:torus-estimate} holds for $q=2$ by the energy estimate.

To obtain the $L^q$-estimate, we apply Lemma \ref{lem:physical-cz}. In detail, given a reference point $Y\in \mathbb{T}^N$, set
\begin{align*}
  W_Y(x, t)=e^{i\mu t}w(Y+\lambda x),\qquad (x, t)\in\mathbb{R}^{d+1}.
\end{align*}
Then
\begin{align}\label{express_nablaWY}
  \nabla_{x, t}W_Y(x, t)=
  \begin{pmatrix}
    \nabla_\lambda w(Y+\lambda x)e^{i\mu t}\\
    i\mu w(Y+\lambda x)e^{i\mu t}
  \end{pmatrix}. 
\end{align}
According to equation \eqref{eq:torus-equation}, by considering the real and imaginary parts separately, $W_Y$ solves
\begin{equation}\label{eq_WY}
  -\diver_{x, t}(\widetilde{A}_Y\nabla_{x, t}W_Y)=\diver_{x, t}\widetilde{H}_Y\quad \text{in }\mathbb{R}^{d+1}
\end{equation}
in the sense of distributions, where
\begin{align*}
  \widetilde{A}_Y(x, t)=
  \begin{pmatrix}
    A(Y+\lambda x) & 0\\
    0 & 1
  \end{pmatrix},\quad \widetilde{H}_Y(x, t)=
  \begin{pmatrix}
    e^{i\mu t}H(Y+\lambda x)\\ 0
  \end{pmatrix}.
\end{align*}
The lifted coefficient matrix $\widetilde{A}_Y$ is also quasiperiodic, that is,
\begin{align*}
  \widetilde{A}_Y(x, t)=B\big((\lambda, 0)\cdot (x, t)^T\big), 
\end{align*}
where $(\lambda, 0)\in \mathbb{R}^{N\times(d+1)}$ is the frequency matrix and $B(Z)=
\begin{pmatrix}
  A(Y+Z) & 0\\
  0 & 1
\end{pmatrix}
$ satisfies the ellipticity, periodicity and H\"{o}lder assumptions as \eqref{eq:ellipticity}--\eqref{eq:holder}. Thus, we can apply Lemma \ref{lem:physical-cz} to obtain for $q>2$
\begin{align*}
  &\quad\pnorm{\nabla_{x, t}W_Y}{q}{B_R\times (-R, R)}\\&\leq C_q\Big\{\pnorm{\nabla_{x, t}W_Y}{2}{B_{5R}\times (-5R, 5R)}+\pnorm{\widetilde{H}_Y}{q}{B_{5R}\times(-5R, 5R)}\Big\}\\&\leq C_q\Big\{\pnorm{\widetilde{H}_Y}{q}{B_{10R}\times(-10R, 10R)}+R^{-1}\anorm{W_Y}{B_{10R}\times(-10R, 10R)}\Big\},
\end{align*}
where the second inequality follows from the Caccioppoli inequality to equation \eqref{eq_WY}. By \eqref{express_nablaWY}, this implies
\begin{align*}
  &\fint_{B_R}|\nabla_\lambda w(Y+\lambda x)|^q+\mu^q|w(Y+\lambda x)|^q \dd x
  \leq C_q\fint_{B_{10R}}|H(Y+\lambda x)|^q \dd x\\
  &\qquad\qquad\qquad\qquad+C_qR^{-q}\bigg(\fint_{B_{10R}}|w(Y+\lambda x)|^2 \dd x\bigg)^{q/2}.
\end{align*}
Integrating this inequality over $Y\in \mathbb{T}^N$ and using Fubini's theorem, one can get
\begin{align}
  \|\nabla_\lambda w\|_q^q+\mu^q\|w\|_q^q
  &\leq C_q\|H\|_q^q+C_q R^{-q}\int_{\mathbb{T}^N}\bigg(\fint_{B_{10R}}|w(Y+\lambda x)|^2 \dd x\bigg)^{q/2}\dd Y\nonumber\\
  &\leq C_q\|H\|_q^q+C_q R^{-q}\|w\|_q^q,
  \label{es_torus-estimate_R}
\end{align}
where H\"{o}lder's inequality is used in the last step.

If we have known $w\in L^q(\mathbb{T}^N)$, by letting $R\rightarrow \infty$ in \eqref{es_torus-estimate_R}, \eqref{eq:torus-estimate} follows immediately. This is true in the case $A$ and $H$ are smooth. In fact, for multi-index $\beta\in \mathbb{N}^N$, the function $v=\partial^\beta w$ satisfies
\begin{align*}
  \mu^2v-\nabla_\lambda\cdot(A\nabla_\lambda v)=\nabla_\lambda\cdot \bigg[\partial^\beta H+\sum_{0<\gamma\leq \beta}
  \begin{pmatrix}
    \beta\\ \gamma
  \end{pmatrix}\partial^\gamma A\nabla_\lambda \partial^{\beta-\gamma} w\bigg]
\end{align*}
in the sense of distributions. Using an induction argument on $\beta$ as well as the energy estimate, one can actually show that $w$ is smooth. For the general case without smoothness assumptions, estimate \eqref{eq:torus-estimate} can be justified by an approximation argument on $A$ and $H$. We leave the details to readers.
\end{proof}

\begin{lemma}\label{lem:phase-regularity}
Let $\X_\mu$ be given by equation \eqref{eq:cell-equation}. For $2<q<\infty$ and $\ell\ge0$, we have
\begin{equation}\label{eq:phase-block-estimate}
\mu\|\Delta_\ell\X_\mu\|_{q}
+\|\nlam\Delta_\ell\X_\mu\|_{q}
+\|\Delta_\ell[A(I+\nlam\X_\mu)]\|_{q}
\le C_q2^{-\alpha\ell}.
\end{equation}
For any $\kappa$ satisfying $0<\kappa<\alpha$, we also have
\begin{equation}\label{eq:phase-block-sup}
\mu\|\Delta_\ell\X_\mu\|_{\infty}
+\|\nlam\Delta_\ell\X_\mu\|_{\infty}
+\|\Delta_\ell[A(I+\nlam\X_\mu)]\|_{\infty}
\le C2^{-\kappa\ell}.
\end{equation}
The constant $C_q$ in \eqref{eq:phase-block-estimate} depends only on $q, \vartheta$ and the character of $A$, while the constant $C$ in \eqref{eq:phase-block-sup} may also depend on $\kappa$.
\end{lemma}
\begin{proof}
  We first apply Lemma \ref{lem:torus-estimate} to equation \eqref{eq:cell-equation} with $H=Ae_j$ to obtain
\begin{equation}\label{eq:cell-Lq}
\mu\|\X_\mu\|_{q}+\|\nabla_\lambda\X_\mu\|_{q}+\|A(I+\nabla_\lambda\X_\mu)\|_{q}\le C_q.
\end{equation}
Recall that $\tau_ZV(Y)=V(Y+Z)$ and $D_ZV=\tau_ZV-V$ for $Z\in\R^N$. Applying $D_Z$ to equation \eqref{eq:cell-equation}, we obtain
\begin{equation}\label{eq_DZX}
  (\mu^2-\nlam\cdot A\nlam)D_Z\X_{\mu,j}
  =\nlam\cdot\bigl\{(D_ZA)(e_j+\tau_Z\nabla_\lambda\X_{\mu,j})\bigr\}.
\end{equation}
Note that, for any $Z\in\mathbb{R}^N$, by the H\"older continuity of $A$ and \eqref{eq:cell-Lq}, it holds that
\begin{equation}\label{es_DZA_q}
\|(D_ZA)(e_j+\tau_Z\nabla_\lambda\X_{\mu,j})\|_{q}
\le\|D_ZA\|_{\infty}\bigl(1+\|\nabla_\lambda\X_{\mu}\|_{q}\bigr)\le C_q|Z|^\alpha.
\end{equation}
Thus, applying Lemma \ref{lem:torus-estimate} to equation \eqref{eq_DZX}, we obtain
\begin{equation}\label{eq:cell-phase-difference}
\mu\|D_Z\X_\mu\|_{q}+\|D_Z\nabla_\lambda\X_\mu\|_{q}\le C_q|Z|^\alpha.
\end{equation}
Besides, since
\begin{align*}
  D_Z[A(I+\nabla_\lambda\X_\mu)]=(D_ZA)(I+\tau_Z\nabla_\lambda\X_\mu)+A D_Z\nabla_\lambda\X_\mu,
\end{align*}
\eqref{es_DZA_q} and \eqref{eq:cell-phase-difference} imply that
\begin{equation}\label{eq:flux-phase-difference}
\|D_Z[A(I+\nabla_\lambda\X_\mu)]\|_{q}\le C_q|Z|^\alpha.
\end{equation}

Now we are ready to prove \eqref{eq:phase-block-estimate}. The case $\ell=0$ follows from \eqref{eq:cell-Lq} as well as the $L^q(\T^N)$ boundedness of $\Delta_0$. Noticing \eqref{eq:cell-phase-difference}--\eqref{eq:flux-phase-difference}, the case $\ell\ge1$ follows from Lemma \ref{diff-block} directly.

To show \eqref{eq:phase-block-sup}, we choose $q>2$ large enough such that $N/q\leq \alpha-\kappa$, depending on $N,\alpha,\kappa$. By Lemma \ref{Bern}, we conclude that
\begin{align*}
&\quad\mu\|\Delta_\ell\X_\mu\|_{\infty}
+\|\Delta_\ell\nabla_\lambda\X_\mu\|_{\infty}
+\|\Delta_\ell[A(I+\nabla_\lambda\X_\mu)]\|_{\infty}\\
&\le C_q2^{N\ell/q}
\bigl\{\mu\|\Delta_\ell\X_\mu\|_{q}
+\|\Delta_\ell\nabla_\lambda\X_\mu\|_{q}
+\|\Delta_\ell[A(I+\nabla_\lambda\X_\mu)]\|_{q}\bigr\}\\
&\le C_q2^{-(\alpha-N/q)\ell}\\&\le C_q2^{-\kappa\ell},
\end{align*}
where \eqref{eq:phase-block-estimate} is used in the second inequality. 
\end{proof}

\subsection{Flux correctors}

Denote
\begin{equation}\label{def_GHmu}
\mathcal G_\mu:=\nlam\X_\mu\quad\text{and}\quad
\mathcal H_\mu:=A(I+\mathcal \nlam\X_\mu).
\end{equation}
For $\mu>0$ and $1\le k,i,j\le d$, we define the flux corrector $\Phi_{\mu,kij}$ as the weak solution of the cell problem
\begin{equation}\label{eq:mu-auxiliary}
(\mu^2-\Delta_\lambda)\Phi_{\mu,kij}
=\nabla_{\lambda,i}\mathcal H_{\mu,kj}
-\nabla_{\lambda,k}\mathcal H_{\mu,ij}\quad \text{in }\T^N,
\end{equation}
where the directional Laplacian is given by
\begin{align}\label{eq:def-directional_laplace}
  \Delta_\lambda:=\nabla_\lambda\cdot \nabla_\lambda=\sum_{\ell=1}^d\nabla_{\lambda,\ell}^2.
\end{align}
Note that by Lemma \ref{lem:torus-estimate} equation \eqref{eq:mu-auxiliary} is well-posed and $\Phi_{\mu,kij}\in H^1_\lambda(\T^N)$ satisfies the energy estimate
\begin{equation*}
\mu\|\Phi_{\mu}\|_2+\|\nlam\Phi_{\mu}\|_2\le C\|\mathcal H_{\mu}\|_2.
\end{equation*}
Besides, $\Phi_\mu$ is skew-symmetric in $k, i$, i.e., $\Phi_{\mu,kij}=-\Phi_{\mu,ikj}$, and its Fourier coefficients are given by
\begin{equation}\label{eq:massive-flux-symbol}
\widehat{\Phi_{\mu,kij}}(n)
=\frac{2\pi i\bigl(
\xi_i\widehat{\mathcal H_{\mu,kj}}(n)
-\xi_k\widehat{\mathcal H_{\mu,ij}}(n)\bigr)}
{\mu^2+4\pi^2|\xi|^2},
\qquad \xi=\xi_\lambda(n),\quad n\in\Z^N.
\end{equation}
Evidently, \eqref{eq:massive-flux-symbol} implies that $\Pi_\lambda \Phi_\mu=0$ by the definition of $\Pi_\lambda$.

Denote
\begin{equation}\label{eq:mu-u-definition}
\U_{\mu,ij}:=\mathcal H_{\mu,ij}
-\sum_{k=1}^d\nabla_{\lambda,k}\Phi_{\mu,kij}\in L^2(\T^N).
\end{equation}
A direct computation shows that
\begin{equation*}
\begin{aligned}
(\mu^2-\Delta_\lambda)\U_{\mu,ij}
&=\mu^2\mathcal H_{\mu,ij}-\Delta_\lambda\mathcal H_{\mu,ij}
-\sum_k\nabla_{\lambda,k}
\bigl(\nabla_{\lambda,i}\mathcal H_{\mu,kj}
-\nabla_{\lambda,k}\mathcal H_{\mu,ij}\bigr)\\
&=\mu^2\mathcal H_{\mu,ij}
-\nabla_{\lambda,i}\sum_k\nabla_{\lambda,k}\mathcal H_{\mu,kj}\\
&=\mu^2(\mathcal H_{\mu,ij}-\mathcal G_{\mu,ij}),
\end{aligned}
\end{equation*}
where the fact $\nlam\cdot\mathcal H_\mu=\mu^2\X_\mu$ from equation \eqref{eq:cell-equation} is used. Combining the above identities with (4.15), we conclude
\begin{equation}\label{eq:resolvent-cell-flux}
\begin{gathered}
(\mu^2-\Delta_\lambda)\X_\mu
=\nlam\cdot(\mathcal H_\mu-\mathcal G_\mu),\\
(\mu^2-\Delta_\lambda)\Phi_{\mu,kij}
=\nabla_{\lambda,i}\mathcal H_{\mu,kj}
-\nabla_{\lambda,k}\mathcal H_{\mu,ij},\\
(\mu^2-\Delta_\lambda)\U_\mu
=\mu^2(\mathcal H_\mu-\mathcal G_\mu).
\end{gathered}
\end{equation}
As in \eqref{eq:massive-flux-symbol}, one can deduce the Fourier coefficients of $\X_\mu$ from the first equation in \eqref{eq:resolvent-cell-flux} and show that $\Pi_\lambda\X_\mu=0$. Besides, by (iii) in Lemma \ref{lem:orbit-averages}, we also have
\begin{align*}
  \Pi_\lambda[\U_\mu(\,\cdot+\lambda x)-\U_\mu]=\Pi_\lambda \U_\mu-\Pi_\lambda \U_\mu=0\quad \text{for any }x\in\mathbb{R}^d.
\end{align*}

\begin{lemma}\label{lem:band-resolvent}
For $\mu>0$ and $\ell\ge0$, we have
\begin{equation}\label{eq:flux-lp-block-sup}
\mu\|\Delta_\ell\Phi_\mu\|_\infty
+\|\nlam\Delta_\ell\Phi_\mu\|_\infty
+\|\Delta_\ell\U_\mu\|_\infty
\le C2^{-\kappa\ell},
\end{equation}
where $C$ depends only on $\kappa, \vartheta$ and the character of $A$.
\end{lemma}

\begin{proof}
For $2<q<\infty$ and $1\le k,i,j\le d$, let
\begin{equation*}
w=\Delta_\ell\Phi_{\mu,kij}\quad\text{and}\quad
H=e_i\Delta_\ell\mathcal H_{\mu,kj}
-e_k\Delta_\ell\mathcal H_{\mu,ij},
\end{equation*}
where $e_i,e_k$ are the coordinate unit vectors in $\R^d$.
Since $\Delta_\ell$ commutes with $\nlam$, equation \eqref{eq:mu-auxiliary} yields
\begin{equation*}
\mu^2w-\nlam\cdot(I\nlam w)=\nlam\cdot H\quad \text{in }\T^N.
\end{equation*}
Applying Lemma~\ref{lem:torus-estimate} with $A=I$ and using \eqref{eq:phase-block-estimate}, we obtain
\begin{equation*}
\mu\|\Delta_\ell\Phi_\mu\|_q
+\|\nlam\Delta_\ell\Phi_\mu\|_q
\le C_q\|\Delta_\ell\mathcal H_\mu\|_q
\le C_q2^{-\alpha\ell}.
\end{equation*}
By \eqref{eq:mu-u-definition}, this yields
\begin{equation*}
\|\Delta_\ell\U_\mu\|_q
\le\|\Delta_\ell\mathcal H_\mu\|_q
+C_d\|\nlam\Delta_\ell\Phi_\mu\|_q
\le C_q2^{-\alpha\ell}.
\end{equation*}
Choose $q>\max\{2,N/(\alpha-\kappa)\}$. Lemma \ref{Bern} implies 
\begin{equation*}
\mu\|\Delta_\ell\Phi_\mu\|_\infty
+\|\nlam\Delta_\ell\Phi_\mu\|_\infty
+\|\Delta_\ell\U_\mu\|_\infty
\le C_q2^{N\ell/q}2^{-\alpha\ell}
\le C2^{-\kappa\ell}.
\end{equation*}
The proof is complete.
\end{proof}

\begin{corollary}\label{err1_coro-bound}
  For $\mu>0$, $\X_\mu, \nlam \X_\mu, \Phi_\mu, \nlam\Phi_\mu, \U_\mu$ are all continuous functions on $\T^N$, and there exists $C_0$, depending only on $\kappa, \vartheta$ and the character of $A$, such that
  \begin{align*}
    \mu\|\X_\mu\|_{\infty}+\|\nlam \X_\mu\|_{\infty}+\mu\|\Phi_\mu\|_\infty+\|\nlam\Phi_\mu\|_\infty+\|\U_\mu\|_\infty\leq C_0.
  \end{align*}
\end{corollary}
\begin{proof}
  Let $V$ be either one of $\mu\X_\mu, \nlam \X_\mu, \mu\Phi_\mu, \nlam\Phi_\mu, \U_\mu$. The series $\sum_{\ell=0}^\infty\Delta_\ell V$ converges uniformly according to \eqref{eq:phase-block-sup} and \eqref{eq:flux-lp-block-sup}. Since each block $\Delta_\ell V$ is smooth and $\sum_{\ell=0}^\infty \Delta_\ell = \mathrm{Id}$, we know $V$ is continuous on $\T^N$ and
  \begin{align*}
    \|V\|_\infty\leq \sum_{\ell=0}^\infty\|\Delta_\ell V\|_\infty\leq C\sum_{\ell=0}^\infty 2^{-\kappa \ell} \le C_0. 
  \end{align*}
  The proof is complete.
\end{proof}

Corollary \ref{err1_coro-bound} only provides uniform boundedness. To continue, we need the following precise smallness estimates in $\mu$ based on double decompositions. 

\begin{lemma}\label{lem:cell-flux-bounds}
  Let $0<\kappa<\alpha$ and $\omega_\lambda=\omega_{\lambda,\kappa}$ be given by \eqref{eq:def-modulus}. For every $\mu>0$ and $r>0$, we have
\begin{equation}\label{eq:phase-smallness}
\begin{aligned}
\mu\bigl(\|\X_\mu\|_{\infty}
+\|\Phi_\mu\|_{\infty}\bigr)
&\le C\omega_\lambda(\mu^{-1}),\\
\sup_{|h|\le r}
\|\U_\mu(\,\cdot+\lambda h)-\U_\mu\|_{\infty}
&\le C\mu r\,\omega_\lambda(\mu^{-1}),
\end{aligned}
\end{equation}
where $C$ depends only on $\kappa, \vartheta$ and the character of $A$.
\end{lemma}

\begin{proof}
  We first note that $\omega_\lambda(t)$ is well-defined. In fact, since $0\le\psi\le1$, it follows from Corollary \ref{cor:counting} that for $t>0$
\begin{equation*}
\omega_\lambda(t)\le\sum_{\ell\ge0}2^{-\kappa\ell}\#\cJ_\ell^\lambda\le C_N\sum_{\ell\ge0}(\ell+1)2^{-\kappa\ell}<\infty.
\end{equation*}
The uniform convergence also implies that $\omega_\lambda(t)$ is a nonnegative continuous function. 

Next we deduce \eqref{eq:phase-smallness} from the claim that, for $V\in C(\T^N)$, $\mu>0$ and $s\in\Z$,
\begin{equation}\label{eq:op-band-resolvent}
\|P_s^\lambda\nlam(\mu^2-\Delta_\lambda)^{-1}V\|_\infty
\le C_d\frac{2^s}{\mu^2+2^{2s}}\|V\|_\infty.
\end{equation}
Due to equations \eqref{eq:resolvent-cell-flux}, for $V=\X_\mu, \Phi_\mu$ or $ \mu^{-2}\nabla_\lambda \mathcal{A}_\mu$, respectively, we have from \eqref{eq:op-band-resolvent} and \eqref{eq:phase-block-sup} that
\begin{align*}
  \|P_s^\lambda\Delta_\ell V\|_\infty\leq \frac{C2^s}{\mu^2+2^{2s}} ( \|\Delta_\ell \cH_\mu\|_\infty + \|\Delta_\ell \G_\mu\|_\infty ) \le C2^{-\kappa\ell}\frac{2^s}{\mu^2+2^{2s}}.
\end{align*}
Note that
\begin{equation*}
\frac{\mu2^s}{\mu^2+2^{2s}}\le\min\{2^s\mu^{-1},2^{-s}\mu\}=\psi(2^s\mu^{-1}).
\end{equation*}
The inequality above yields
\begin{equation}\label{eq:band-cell-flux}
\begin{aligned}
\mu\bigl(\|P_s^\lambda\Delta_\ell\X_\mu\|_\infty
+\|P_s^\lambda\Delta_\ell\Phi_\mu\|_\infty\bigr)
&\le C2^{-\kappa\ell}\psi(2^s\mu^{-1}),\\
\|\nlam P_s^\lambda\Delta_\ell\U_\mu\|_\infty
&\le C\mu\,2^{-\kappa\ell}\psi(2^s\mu^{-1}).
\end{aligned}
\end{equation}
Moreover, since
\begin{equation*}
\Pi_\lambda\Phi_\mu=\Pi_\lambda\X_\mu=\Pi_\lambda[\U_\mu(\,\cdot+\lambda h)-\U_\mu]=0
\end{equation*}
for any $h\in\R^d$, \eqref{eq:lp-reconstruction} implies
\begin{equation}\label{eq:potential-band-reconstruction}
 V= \sum_{\ell \ge 0} \sum_{s\in\cJ_\ell^\lambda}P_s^\lambda\Delta_\ell V,
\quad\text{for }
V=\X_\mu,\ \Phi_\mu,\ \U_\mu(\,\cdot+\lambda h)-\U_\mu.
\end{equation}
Consequently, it follows from \eqref{eq:band-cell-flux}--\eqref{eq:potential-band-reconstruction} that
\begin{equation*}
\begin{aligned}
\mu\bigl(\|\X_\mu\|_\infty+\|\Phi_\mu\|_\infty\bigr)
&\le\mu\sum_{\ell\ge0}\sum_{s\in\cJ_\ell^\lambda}
\bigl(\|P_s^\lambda\Delta_\ell\X_\mu\|_\infty
+\|P_s^\lambda\Delta_\ell\Phi_\mu\|_\infty\bigr)\\
&\le C\sum_{\ell\ge0}2^{-\kappa\ell}
\sum_{s\in\cJ_\ell^\lambda}\psi(2^s\mu^{-1})
=C\omega_\lambda(\mu^{-1}),
\end{aligned}
\end{equation*}
and
\begin{equation*}
\begin{aligned}
\sup_{|h|\le r}
\|\U_\mu(\,\cdot+\lambda h)-\U_\mu\|_\infty
&\le\sum_{\ell\ge0}\sum_{s\in\cJ_\ell^\lambda}\sup_{|h|\le r}\|P_s^\lambda\Delta_\ell
[\U_\mu(\,\cdot+\lambda h)-\U_\mu]\|_\infty\\&\le \sum_{\ell\ge0}\sum_{s\in\cJ_\ell^\lambda}r\|\nabla_\lambda P_s^\lambda \Delta_\ell\mathcal{A}_\mu\|_\infty\\
&\le C\mu r\sum_{\ell\ge0}2^{-\kappa\ell}
\sum_{s\in\cJ_\ell^\lambda}\psi(2^s\mu^{-1})
=C\mu r\,\omega_\lambda(\mu^{-1}),
\end{aligned}
\end{equation*}
where we have used the mean value theorem in the second inequality.

It remains to prove \eqref{eq:op-band-resolvent}. Set 
$$\widetilde{\mathcal{K}} =\mathcal F^{-1}\left(\frac{\rho(\varsigma)}{\mu^2+4\pi^2 2^{2s}|\varsigma|^2}\right),$$ 
where $\mathcal{F}$ denotes the Fourier transform in $\R^d$. Note that $\widetilde{\mathcal{K}}$ is a Schwartz function. Then, calculating in the sense of distributions,
\begin{align}\label{iden:convolution}
\begin{split}
    (\mu^2-\Delta_\lambda)^{-1}P_s^\lambda V(Y)&=\sum_{n\in\mathbb{Z}^N}\widehat{V}(n)e^{2\pi in\cdot Y}\frac{\rho(2^{-s}\xi_\lambda(n))}{\mu^2+4\pi^2|\xi_\lambda(n)|^2}\\
    &=\sum_{n\in\mathbb{Z}^N}\widehat{V}(n)e^{2\pi in\cdot Y}\int_{\mathbb{R}^d}\widetilde{\mathcal{K}}(y)e^{-2\pi i y\cdot \lambda^Tn/2^s}\dd y\\
    &=\int_{\mathbb{R}^d}\widetilde{\mathcal{K}}(y)\sum_{n\in\mathbb{Z}^N}\widehat{V}(n)e^{2\pi in\cdot(Y-\lambda y/2^s)}\dd y\\
    &=\int_{\mathbb{R}^d}\widetilde{\mathcal{K}}(y)V(Y-\lambda y/2^s)\dd y.
\end{split}
\end{align}
Therefore, we obtain from integration by parts
\begin{align*}
 \nabla_\lambda(\mu^2-\Delta_\lambda)^{-1} P^\lambda_sV(Y)&=-2^{s}\int_{\mathbb R^d} \widetilde{\mathcal{K}}(y)\nabla_y[V(Y-\lambda y/2^s)]\dd y\\
 &=2^s\int_{\mathbb R^d}\nabla \widetilde{\mathcal{K}}(y)V(Y-\lambda y/2^s)\dd y,
\end{align*}
which, by Minkowski's inequality, yields
\begin{equation*}
 \|\nabla_\lambda(\mu^2-\Delta_\lambda)^{-1} P^\lambda_sV(Y)\|_\infty
 \le 2^s\|\nabla \widetilde{\mathcal{K}}\|_{L^1(\mathbb R^d)}\|V\|_\infty.
\end{equation*}
Note that using the support of $\rho$
\begin{align*}
 \|\nabla \widetilde{\mathcal{K}}\|_{L^1(\R^d)}&=\|(1+|y|^2)^{-d}\cdot(1+|y|^2)^d\nabla \widetilde{\mathcal{K}}\|_{L^1(\R^d)}\le C_d\|(1+|y|^2)^d\nabla \widetilde{\mathcal{K}}\|_{L^\infty(\R^d)}\\
 &\le C_d \left\|(\mathrm{Id}-\Delta_\varsigma)^d
 \left(\frac{\varsigma\rho(\varsigma)}{\mu^2+4\pi^22^{2s}|\varsigma|^2}\right)
 \right\|_{L^1(\R^d)}\le \frac{C_d}{\mu^2+2^{2s}},
\end{align*}
where $C_d$ depends only on $d$. This proves the claim \eqref{eq:op-band-resolvent}. 
\end{proof}

\subsection{Low-high frequency separation}
\label{subsec:lowhigh}

The idea of low-high frequency separation can already be seen in the two-scale homogenization problem with coefficient matrix in the form of $A(x,x/\eps)$ and $A(x,y)$ being periodic in the fast variable $y$. Let $\chi(x,y)$ and $\phi(x,y)$ be the correctors and flux correctors constructed in \cite{NSX20}, then
\begin{equation}\label{eq:simple-lowhigh}
    A(I +\nabla_y \chi)(x,x/\eps) = \big[ \widehat{A}(x) - \eps (\diver_x \phi)(x,x/\eps) \big] + \diver (\eps \phi(x,x/\eps)).
\end{equation}
If $A(x,y)$ is smooth, then the terms in brackets above are Lipschitz independent of $\eps$, which can be interpreted as a low frequency part. The divergence term in \eqref{eq:simple-lowhigh} oscillates rapidly at $\eps$ scale, enclosing high frequencies. This is a simple form of low-high frequency separation. In the following, we apply a similar idea to general quasiperiodic problems.

Restricting $\X_\mu(Y)$ and $\Phi_\mu(Y)$ onto the physical orbit $x \mapsto \lambda x$, we recover the correctors and flux correctors in $\R^d$,
\begin{equation}\label{eq:physical-fields}
\chi_\mu(x):=\X_\mu(\lambda x),\qquad
\phi_\mu(x):=\Phi_\mu(\lambda x),\qquad x\in\R^d.
\end{equation}
By \eqref{eq:cell-equation}
\begin{align}
  \label{eq_chi_mu}
  \diver(A_\lambda(I+\nabla \chi_\mu))=\mu^2\chi_\mu,
\end{align}
and by \eqref{eq:mu-u-definition}
\begin{equation}\label{eq:primal-flux}
A_\lambda(I+\nabla\chi_\mu)=\mathcal A_\mu(\lambda x) + \diver\phi_\mu
\quad\text{in }\R^d,
\end{equation}
where \begin{equation*}
(\nabla\chi_\mu)_{ij}=\partial_i\chi_{\mu,j},
\qquad
(\diver\phi_\mu)_{ij}=\sum_{k=1}^d\partial_k\phi_{\mu,kij}.
\end{equation*}
Moreover, $\phi_{\mu,kij}$ is skew-symmetric in the first two indices, i.e., $\phi_{\mu,kij}=-\phi_{\mu,ikj}$.

Now we would like to show \eqref{eq:primal-flux} is a low-high frequency separation in the same spirit as in two-scale problem. To see this, we need to show that $\mathcal A_\mu$ oscillates mildly below a certain scale.
Fix $x_0 \in \R^d$ and let
\begin{equation}\label{def_overlineA}
\overline{\mathcal A}:=\U_\mu(\lambda x_0),\qquad \mathcal R(x):=\U_\mu(\lambda x)-\overline{\mathcal A},\qquad x\in\R^d.
\end{equation}
Then \eqref{eq:primal-flux} can be rewritten as
\begin{equation}\label{eq:primal-comparison-flux}
A_\lambda(I+\nabla\chi_\mu)=\overline{\mathcal A}+\mathcal R+\diver\phi_\mu
\quad\text{in }\R^d,
\end{equation}
We have from Corollary \ref{err1_coro-bound} and Lemma \ref{lem:cell-flux-bounds} the following estimates. 
\begin{corollary}\label{coro:physical}
Let $r>0$ and $\mu=r^{-1}$. Then $\chi_\mu, \phi_\mu\in C^{1}_b(\R^d)$ and $\mathcal R\in C_b(\R^d;\R^{d\times d})$. 
Moreover, there exists a constant $C_0$, depending only on $\kappa, \vartheta$ and the character of $A$, such that
\begin{equation}\label{eq:decomposition-uniform}
\|\overline{\mathcal A}\|+\|\nabla\chi_\mu\|_{L^\infty(\R^d)}+\|\nabla\phi_\mu\|_{L^\infty(\R^d)}\le C_0,
\end{equation}
and
\begin{equation}\label{eq:decomposition-small}
\|\mathcal R\|_{L^\infty(B_{4r}(x_0))}
+r^{-1}\bigl(\|\chi_\mu\|_{L^\infty(B_{4r}(x_0))}
+\|\phi_\mu\|_{L^\infty(B_{4r}(x_0))}\bigr)
\le C_0\omega_\lambda(r),
\end{equation}
where $\|\overline{\mathcal A}\|$ denotes the operator norm of the matrix $\overline{\mathcal A}$. 
\end{corollary}

From \eqref{eq:decomposition-small} we see that the remainder term $\mathcal{R}(x)$ is small, if $\omega_\lambda(r)$ is small and $x\in B_{4r}(x_0)$. This explains the low-high frequency separation and explains why we can also approximate the quasiperiodic operator in most scales as $r$ varies by elliptic operators with constant coefficients.

The following proposition shows that the constant matrix $\overline{\mathcal{A}}$ possesses the ellipticity if $\omega_\lambda(r)$ is small. 
\begin{proposition}\label{prop:mu-decomposition}
Let $r>0$ and $\mu=r^{-1}$. There exists a constant $C_0$, depending only on $\kappa, \vartheta$ and the character of $A$, such that
\begin{equation}\label{eq:almost-ellipticity}
\zeta\cdot\overline{\mathcal A}\zeta
\ge\bigl(\Lambda^{-1}-C_0\omega_\lambda(r)\bigr)|\zeta|^2
\qquad\text{for any }\zeta\in\R^d.
\end{equation}
In particular, $\overline{\mathcal A}$ is uniformly elliptic if $C_0\omega_\lambda(r)\le(2\Lambda)^{-1}$.
\end{proposition}

\begin{proof}
    We first prove the ellipticity of $  \A_\mu(\lambda x)$ in the average sense. Let $\eta\in C_c^\infty(B_1)$ be a nonnegative cutoff function such that $\int_{B_1}\eta=1$ and set $\eta_r(x)=r^{-d}\eta((x-x_0)/r)$. Then, $\eta_r$ is supported in $B_r$ and
\begin{equation*}
\int_{B_{r}(x_0)}\eta_r=1,
\qquad
\|\nabla\eta_r\|_{L^1(B_{r}(x_0))}\le Cr^{-1}.
\end{equation*}
Let $\zeta \in \R^d$. By \eqref{eq:primal-flux} and integration by parts
\begin{equation*}
\begin{aligned}
    \int \zeta \cdot \mathcal{A}_\mu(\lambda\cdot) \zeta\eta_r & = \int \zeta\cdot A_\lambda(I + \nabla \chi_\mu) \zeta\eta_r  - \int \zeta\cdot \diver(\phi_\mu) \zeta\eta_r \\
    & = \int (I + \nabla \chi_\mu)\zeta\cdot A_\lambda(I + \nabla \chi_\mu) \zeta\eta_r + \int \mu^2 (\chi_\mu \cdot \zeta)^2 \eta_r \\
    & \qquad + \int  (\chi_\mu \cdot \zeta) \nabla \eta_r \cdot A_\lambda(I + \nabla \chi_\mu) \zeta + \int \phi_\mu:(\nabla \eta_r \otimes \zeta\otimes \zeta) \\
    &=: J_1 + J_2 + J_3 + J_4,
\end{aligned}
\end{equation*}
where we have used \eqref{eq_chi_mu} in the last equation. For $J_1$, using the ellipticity of $A_\lambda$,
\begin{equation*}
\begin{aligned}
    J_1 & \ge \Lambda^{-1} \int |(I + \nabla \chi_\mu)\zeta|^2 \eta_r \\
    & = \Lambda^{-1} \int (|\zeta|^2 + |\nabla \chi_\mu\zeta|^2) \eta_r  - 2\Lambda^{-1} \int (\chi_\mu \cdot \zeta) (\nabla \eta_r \cdot \zeta) \\
    & \ge \Lambda^{-1} |\zeta|^2 - C\Lambda^{-1} |\zeta|^2 r^{-1} \| \chi_\mu \|_{L^\infty(\R^d)} \\
    & \ge \Lambda^{-1} |\zeta|^2 - C\Lambda^{-1} |\zeta|^2 \omega_\lambda(r),
\end{aligned}
\end{equation*}
where we have used \eqref{eq:decomposition-small} in the last inequality. We simply drop $J_2$ since it is nonnegative. For $J_3$, we use \eqref{eq:decomposition-uniform} and \eqref{eq:decomposition-small} to obtain
\begin{equation*}
    |J_3| \le C|\zeta|^2 \| \chi_\mu \|_{L^\infty(\R^d)} (1+\| \nabla \chi_\mu \|_{L^\infty(\R^d)}) \| \nabla \eta_r \|_{L^1(B_r(x_0))} \le C|\zeta|^2 \omega_\lambda(r).
\end{equation*}
For $J_4$, we use \eqref{eq:decomposition-small} to get
\begin{equation*}
    |J_4| \le C|\zeta|^2 \| \phi_\mu \|_{L^\infty(\R^d)} \| \nabla \eta_r \|_{L^1(\R^d)} \le C|\zeta|^2 \omega_\lambda(r).
\end{equation*}
Combining the above estimates leads to the averaged ellipticity
\begin{equation*}
    \int \zeta \cdot \mathcal{A}_\mu(\lambda\cdot) \zeta\eta_r \ge |\zeta|^2 ( \Lambda^{-1} - C \omega_\lambda(r) ).
\end{equation*}

Now we show that the averaged ellipticity of $\A_\mu(\lambda x)$ implies the pointwise ellipticity. For any given $x_0 \in \R^d$, by \eqref{def_overlineA} and \eqref{eq:decomposition-small}
\begin{equation*}
\begin{aligned}
    \zeta \cdot \overline{\A} \zeta & = \int  \zeta \cdot \overline{\A} \zeta \eta_r = \int \zeta \cdot \mathcal{A}_\mu(\lambda x) \zeta \eta_r - \int \zeta \cdot \mathcal{R} \zeta  \eta_r \\
    & \ge |\zeta|^2 ( \Lambda^{-1} - C \omega_\lambda(r) ) - |\zeta|^2 \| \mathcal{R}\|_{L^\infty(B_r(x_0))} \\
    & \ge |\zeta|^2 ( \Lambda^{-1} - C \omega_\lambda(r) ).
\end{aligned}
\end{equation*}
This completes the proof.
\end{proof}

\section{Uniform Lipschitz estimates}
\label{sec.proveLip}

In this section, we prove Theorem \ref{thm:main} by applying Proposition \ref{prop:criterion}. The next proposition establishes the comparison estimate in Proposition \ref{prop:criterion} via a kind of duality method. 

\begin{proposition}\label{prop:all-radii}
Let $\lambda\in\R^{N\times d}$, $r>0$, and $x_0\in\R^d$. Suppose that $u\in H^1(B_{r}(x_0))$ is $A_\lambda$-harmonic. Then there exist a uniformly elliptic constant matrix $\mathcal{M}$ and an $\mathcal{M}$-harmonic function $v\in L^2(B_r(x_0))$ such that
\begin{equation}\label{eq:all-radii-comparison}
r^{-1}\anorm{u-v}{B_r(x_0)}
\le C\omega_\lambda(r)\anorm{\nabla u}{B_{r}(x_0)}.
\end{equation}
The ellipticity constant of $\mathcal{M}$ and the constant $C$ in \eqref{eq:all-radii-comparison} depend only on $\kappa, \vartheta$ and the character of $A$.
\end{proposition}
\begin{proof}
Let $\mu=r^{-1}$. Choose $0<\varepsilon_0<1/4$ so that $C_0\varepsilon_0\le(2\Lambda)^{-1}$, where $C_0$ is the constant given in Proposition \ref{prop:mu-decomposition}. If $\omega_\lambda(r)>\varepsilon_0$, let $\mathcal{M}=I$ and $v=\fint_{B_r(x_0)}u$. By Poincar\'e's inequality, we have
\begin{equation*}
r^{-1}\anorm{u-v}{B_r(x_0)}
\le C\anorm{\nabla u}{B_{r}(x_0)}
\le C\varepsilon_0^{-1}\omega_\lambda(r)\anorm{\nabla u}{B_{r}(x_0)}.
\end{equation*}
Such $r$'s are bad scales at which we do not have a good approximation.

It remains to consider the case $0\leq \omega_\lambda(r)\le\varepsilon_0$. Note that $A^T$ satisfies the same assumptions as $A$, where $A^T$ denotes the transpose of the matrix $A$. We can construct the corresponding quantities of \eqref{eq:physical-fields}--\eqref{def_overlineA} as in Section \ref{sec.corrector} with $A^T$ in place of $A$, denoted by $\chi_\mu^*$, $\phi_\mu^*$, $\overline{\mathcal A}^*$ and $\mathcal R^*$. Then by \eqref{eq:primal-comparison-flux}
\begin{equation}\label{eq:physical-decomposition}
A_\lambda^T(I+\nabla\chi_\mu^*)=\overline{\mathcal A}^*+\mathcal R^*+\diver\phi_\mu^*,
\end{equation}
and the estimates \eqref{eq:decomposition-uniform}--\eqref{eq:almost-ellipticity} hold for $\chi_\mu^*, \phi_\mu^*, \overline{\mathcal A}^*$ and $\mathcal R^*$, with the same $C_0$ and $\omega_\lambda$. Since $C_0\omega_\lambda(r)\le(2\Lambda)^{-1}$, we have
\begin{align*}
(2\Lambda)^{-1}|\zeta|^2\le\zeta\cdot \overline{\mathcal A}^*\zeta\qquad\text{for }\zeta\in\R^d.
\end{align*}
Using \eqref{eq:physical-decomposition}, one can see that
\begin{align*}
  -\diver((\overline{\mathcal{A}}^*)^T \nabla u)=\diver((\mathcal{R}^*)^T\nabla u)+\diver((\diver\phi_\mu^*)^T \nabla u)-\diver([A_\lambda^T (I+\nabla \chi_\mu^*)]^T\nabla u).
\end{align*}
Since $\overline{\mathcal{A}}^*$ is constant, it holds that
\begin{align*}
  -\diver((\overline{\mathcal{A}}^*)^T \nabla u)=-\diver(\overline{\mathcal{A}}^* \nabla u).
\end{align*}
On the other hand, calculating in the sense of distributions, we obtain from the skew-symmetry of $\phi_\mu^*$ and the $A_\lambda$-harmonicity that
\begin{align*}
  \diver((\diver\phi_\mu^*)^T \nabla u) & =\sum_{i,j,k=1}^d\partial_j(\partial_k\phi_{\mu,kij}^* \partial_i u)\\
  & = \sum_{i,j,k=1}^d \partial_{jk}^2[\phi_{\mu,kij}^*\partial_i u]-\sum_{i,j,k=1}^d\partial_j[\phi_{\mu,kij}^*\partial_{ki} u]\\
  &=\sum_{i,j,k=1}^d\partial_{jk}^2[\phi_{\mu,kij}^*\partial_{i} u],
\end{align*}
and
\begin{align*}
  -\diver([A_\lambda^T (I+\nabla \chi_\mu^*)]^T\nabla u)&=-\diver(A_\lambda \nabla u)-\diver((\nabla \chi_\mu^*)^T A_\lambda\nabla u)\\&=-\sum_{j, k=1}^d\partial_{jk}^2[\chi_{\mu,j}^*(A_\lambda\nabla u)_k].
\end{align*}
Therefore, it holds in the sense of distributions
\begin{align*}
  \begin{split}
    -\diver(\overline{\mathcal{A}}^* \nabla u) &= \sum_{j, k=1}^d\partial_k(\mathcal{R}^*_{jk}\partial_j u)+\sum_{i,j,k=1}^d\partial_{jk}^2[\phi_{\mu,kij}^*\partial_{i} u] \\&\quad-\sum_{j, k=1}^d\partial_{jk}^2[\chi_{\mu,j}^*(A_\lambda\nabla u)_k]\quad \text{in } B_r(x_0).
  \end{split}
\end{align*}

Now, let $w_1$ be the solutions of
\begin{align*}
  -\diver(\overline{\mathcal{A}}^* \nabla w_1) = \sum_{j, k=1}^d\partial_k(\mathcal{R}^*_{jk}\partial_j u \mathbf{1}_{B_r(x_0)})\quad \text{in }\R^d
\end{align*}
such that $\int_{B_r(x_0)}w_1=0$, and $w_2$ the solution of
\begin{align*}
-\diver(\overline{\mathcal{A}}^* \nabla w_2)&=\sum_{j,k=1}^d\partial_{jk}^2\bigg[\sum_{i=1}^d\phi_{\mu,kij}^*\partial_{i} u\mathbf{1}_{B_r(x_0)}-\chi_{\mu,j}^*(A_\lambda\nabla u)_k\mathbf{1}_{B_r(x_0)}\bigg]\quad \text{in }\R^d
\end{align*}
given by the Newtonian potential. By the energy estimate and Poincar\'{e}'s inequality,
\begin{align*}
  \|w_1\|_{L^2(B_r(x_0))}&\leq Cr\|\nabla w_1\|_{L^2(B_r(x_0))}\leq C r\|\mathcal{R}^*\|_{L^\infty(B_r(x_0))}\|\nabla u\|_{L^2(B_r(x_0))}\\&\leq Cr\omega_\lambda(r)\|\nabla u\|_{L^2(B_{r}(x_0))},
\end{align*}
and, by the Calder\'{o}n-Zygmund estimate,
\begin{align*}
  \|w_2\|_{L^2(B_r(x_0))}&\leq C \bigl(\|\phi_\mu^*\|_{L^\infty(B_r(x_0))}+\|\chi_\mu^*\|_{L^\infty(B_r(x_0))}\bigr)\|\nabla u\|_{L^2(B_r(x_0))}\\&\leq Cr\omega_\lambda(r)\|\nabla u\|_{L^2(B_{r}(x_0))},
\end{align*}
where we have used the estimate \eqref{eq:decomposition-small} for $\chi_\mu^*, \phi_\mu^*, \mathcal{R}^*$. Set $v=u-w_1-w_2\in L^2(B_{r}(x_0))$. Since $\overline{\mathcal A}^*$ is constant, it follows that $$-\diver(\overline{\mathcal{A}}^*\nabla v)=0 \quad\text{in }B_{r}(x_0).$$ This implies that $v$ is smooth in $B_{r}(x_0)$ and, by the estimates of $w_1, w_2$,
\begin{align*}
\|u-v\|_{L^2(B_{r}(x_0))}=\|w_1+w_2\|_{L^2(B_{r}(x_0))}\le Cr\omega_\lambda(r)\|\nabla u\|_{L^2(B_{r}(x_0))}.
\end{align*}
Dividing by $r|B_r|^{1/2}$, we obtain
\begin{equation*}
  r^{-1}\anorm{u-v}{B_r(x_0)}\le C\omega_\lambda(r)\anorm{\nabla u}{B_{r}(x_0)}.
\end{equation*}
Therefore, by taking $\mathcal{M}=\overline{\mathcal A}^*$ and the preceding $v$, we finish the proof. 
\end{proof}

We now verify the Dini condition for the modulus $\omega_\lambda$ in \eqref{eq:all-radii-comparison}.

\begin{lemma}\label{lem:common-modulus}
Let $\kappa>0$, $\lambda\in \R^{N\times d}$ and $\omega_\lambda=\omega_{\lambda, \kappa}$ be given by \eqref{eq:def-modulus}. Then $\omega_\lambda$ satisfies the Dini condition \eqref{eq:dini-functional} with $D_0$ depending only on $N, \kappa$.
\end{lemma}
\begin{proof}
For $j\in\Z$, the change of variables $t=2^jr$ yields
\begin{equation}\label{eq:single-band-dini}
\int_0^\infty\psi(2^jr)\,\frac{\dd r}{r}
=\int_0^\infty\psi(t)\,\frac{\dd t}{t}
=\int_0^1 \dd t
+\int_1^\infty\frac{\dd t}{t^2}
=2.
\end{equation}
Since $\psi\ge0$, it follows from Fubini's theorem and \eqref{eq:single-band-dini} that
\begin{equation*}
\begin{aligned}
\Dini(\omega_\lambda)
&=\int_0^\infty\sum_{\ell\ge0}2^{-\kappa\ell}
 \sum_{j\in\cJ_\ell^\lambda}\psi(2^jr)\,\frac{\dd r}{r}
 =\sum_{\ell\ge0}2^{-\kappa\ell}
 \sum_{j\in\cJ_\ell^\lambda}
 \int_0^\infty\psi(2^jr)\,\frac{\dd r}{r}\\
&=2\sum_{\ell\ge0}2^{-\kappa\ell}\#\cJ_\ell^\lambda
\le C_N\sum_{\ell\ge0}(\ell+1)2^{-\kappa\ell}=C_{N,\kappa}<\infty,
\end{aligned}
\end{equation*}
where we have also used Corollary \ref{cor:counting}.
\end{proof}

\begin{proof}[Proof of Theorem \ref{thm:main}]
Fix $\lambda\in\R^{N\times d}$, $R>0$, and let $\kappa=\alpha/2$. By \eqref{eq:ellipticity} and Proposition \ref{prop:all-radii}, both $A_\lambda$ and the comparison matrices $\mathcal M$ satisfy \eqref{eq:ellipticity} with the ellipticity constant depending only on $\vartheta$ and the character of $A$. It follows from Proposition \ref{prop:all-radii} and Lemma \ref{lem:common-modulus} that the conditions in Proposition \ref{prop:criterion} are all satisfied with the constants $D_0$ and $C_*$ depending only on $\vartheta$ and the character of $A$. Therefore, for $x_0\in B_{R/2}$, we can apply Proposition \ref{prop:criterion} to $u(x_0+\,\cdot\,)$ on $B_{R/2}$ to get
\begin{equation*}
\begin{aligned}
\sup_{0<t<R/2}\anorm{\nabla u}{B_t(x_0)}
\le C\Bigl\{\anorm{\nabla u}{B_R}
+R\pnorm{f}{p}{B_R}+R^\gamma[g]_{C^{0,\gamma}(B_R)}\Bigr\}.
\end{aligned}
\end{equation*}
Letting $t\rightarrow 0$, we prove \eqref{eq:Lipschitz}.
\end{proof}

We conclude this section with a Liouville theorem.
\begin{corollary}[Zeroth-order Liouville theorem]\label{cor:zero-liouville}
Let $A$ satisfy the same assumptions as in Theorem \ref{thm:main}, and let $\lambda\in\R^{N\times d}$.
Suppose that $u\in H^1_{\mathrm{loc}}(\R^d)$ is a weak solution of
$-\diver(A_\lambda\nabla u)=0$ in $\R^d$ and satisfies
\begin{equation}\label{eq:liouville-sublinear}
\lim_{R\to\infty}\frac{1}{R}\anorm{u}{B_R}=0.
\end{equation}
Then $u$ is constant.
\end{corollary}

\begin{proof}
Fix $\rho>0$. For $R>4\rho$, by Theorem \ref{thm:main} and Caccioppoli's inequality, we have
\begin{equation*}
\|\nabla u\|_{L^\infty(B_\rho)}
\le\|\nabla u\|_{L^\infty(B_{R/4})}
\le C\anorm{\nabla u}{B_{R/2}}
\le\frac{C}{R}\anorm{u}{B_R}.
\end{equation*}
Here $C$ depends only on the character of $A$ and is independent of $R$. Letting $R\to\infty$ and using \eqref{eq:liouville-sublinear}, we obtain $\nabla u=0$ in $B_\rho$. Since $\rho>0$ is arbitrary, $u$ is constant in $\R^d$.
\end{proof}

\section{Quantitative estimates of true correctors}
\label{sec.lineargrowth}

In this section, we will quantify the sublinear growth rates of correctors and flux correctors in terms of $\omega_{\lambda,m+\alpha}$ stated in Theorem \ref{thm:correctors}, which is shown to be optimal by examples. Throughout this section, we assume $A$ satisfies \eqref{eq:ellipticity}, \eqref{eq:periodicity} and \eqref{eq:s6-assumptions}.

\subsection{Sublinear growth}

To construct the true correctors, we take the limit of the approximate correctors in \eqref{eq:cell-equation} as $\mu \to 0$. However, $\X_\mu$ need not converge. Instead, we first consider the convergence of $\nabla_\lambda \X_\mu$. 

\begin{lemma}\label{lem:corrector-existence}
    There exists a continuous field $\mathcal{G} \in C^{m, \kappa}(\T^N;\R^{d\times d})$ for any $\kappa \in (0,\alpha)$, such that, as $\mu \to 0$
    \begin{equation*}
        \mathcal{G}_\mu \to \mathcal{G}, \quad \mathcal{H}_\mu\to \cH:= A(I + \mathcal{G}) \quad \text{uniformly on } \T^N,
      \end{equation*}
      where $\mathcal{G}_\mu, \mathcal{H}_\mu$ are given by \eqref{def_GHmu}.

      Moreover, $(\nabla_\lambda)_i \G_{jk} = (\nabla_\lambda)_j \mathcal{G}_{ik}$, $\Pi_\lambda \mathcal{G} = 0, \nabla_\lambda \cdot \cH = 0$, and for any $2\le q< \infty$, $\ell\geq 0$
    \begin{align}\label{eq:G+H}
        \| \Delta_\ell \mathcal{G} \|_q + \| \Delta_\ell \cH  \|_q \le C_q 2^{-(m+\alpha) \ell}.
    \end{align}
\end{lemma}

\begin{proof}
    Recall $\G_\mu = \nabla_\lambda \X_\mu$ and $\cH_\mu = A(I + \nabla_\lambda \X_\mu)$. The key is to prove that, for any $2\le q<\infty$ and $\ell \ge 0$,
    \begin{equation}\label{eq:higher-block}
    \|\Delta_\ell \G_\mu\|_q + \|\Delta_\ell \cH_\mu\|_q
    \le C_q2^{-(m+\alpha)\ell}.
    \end{equation}

     If \eqref{eq:higher-block} holds, then by Bernstein's inequality in Lemma \ref{Bern}, we have for $k\geq 0$
    \begin{equation*}
      \| \nabla^k \Delta_\ell \G_\mu \|_\infty \le C^{k+1} 2^{k \ell + N \ell /q} \| \Delta_\ell \G_\mu \|_q \le C 2^{(k + N/q -m - \alpha)\ell}.
    \end{equation*}
By choosing $q$ large enough such that $N/q < \alpha$, we know that
\begin{align*}
  \nabla^k\mathcal{G}_\mu=\sum_{\ell \in \N} \nabla^k \Delta_\ell\mathcal{G}_\mu \text{ converges uniformly for } k \le m,
\end{align*}
i.e., $\{\mathcal{G}_\mu\}\subset C^{m}(\T^N;\R^{d\times d})$ with a uniform bound. Further, for each $Z\in \R^N$
  \begin{equation*}
    \begin{aligned}
        \| D_Z \nabla^m \Delta_\ell \G_\mu \|_\infty \le C|Z| \| \nabla^{m+1} \Delta_\ell \G_\mu\|_\infty \le C|Z| 2^{(N/q+1-\alpha) \ell},
    \end{aligned}
  \end{equation*}
  and
  \begin{equation*}
    \| D_Z \nabla^m \Delta_\ell \G_\mu \|_\infty \le 2 \|\nabla^m \Delta_\ell \G_\mu \|_\infty \le C2^{(N/q-\alpha)\ell}.
  \end{equation*}
Letting $\ell_0\in \N$ such that $2^{-\ell_0-1}\leq |Z|\leq 2^{-\ell_0}$, we calculate
    \begin{align*}
        \| D_Z \nabla^m \G_\mu \|_\infty\le C\sum_{\ell\leq \ell_0} 2^{(N/q+1-\alpha)\ell}|Z|+C\sum_{\ell>\ell_0} 2^{(N/q-\alpha)\ell} \le C|Z|^{\alpha - N/q}.
    \end{align*}
    For any $\kappa \in (0,\alpha)$, by choosing $q$ large enough such that $\kappa = \alpha - N/q$, we know $\{\G_\mu\}$ is uniformly bounded in $C^{m,\kappa}(\T^N)$. Hence, by the Arzel\'{a}-Ascoli theorem, we can find a subsequence of $\{ \G_\mu \}$ that converges uniformly to $\G \in C^{m,\kappa}(\T^N)$ as $\mu \to 0$. Up to a subsequence, we also have $\cH_\mu \to \cH= A(I + \G)$ uniformly on $\T^N$. Taking the limit in \eqref{eq:higher-block}, we get \eqref{eq:G+H}, and it is easy to verify that in the sense of distributions
    \begin{align*}
      (\nabla_\lambda)_i \G_{jk} = (\nabla_\lambda)_j \mathcal{G}_{ik},\quad \Pi_\lambda \mathcal{G} = 0,\quad \nabla_\lambda \cdot \cH = 0. 
    \end{align*}
Note that $\nabla_\lambda \cdot \mathcal{H}=0$ implies that for any $v\in C^\infty(\T^N)$
\begin{align*}
  \int_{\T^N} \nabla_\lambda v\cdot A (I+\mathcal{G})=0. 
\end{align*}
If there exists another $\tilde{\mathcal{G}}$ satisfying the same properties, then we have
\begin{align*}
  \int_{\T^N} \nabla_\lambda v\cdot A (\mathcal{G}-\tilde{\mathcal{G}})=0\quad \text{for any }v\in C^\infty(\T^N).
\end{align*}
As $\mathcal{G}-\widetilde{\mathcal{G}}$ is a uniform limit of functions in $\{\nabla_\lambda v: v\in C^\infty(\T^N)\}$, we deduce that
\begin{align*}
  \int_{\T^N} (\mathcal{G}-\tilde{\mathcal{G}})\cdot A (\mathcal{G}-\tilde{\mathcal{G}})=0,
\end{align*}
which, by the ellipticity of $A$, implies that $\mathcal{G}=\tilde{\mathcal{G}}$ in $\T^N$. This proves the uniqueness of the limit $\mathcal{G}$. Together with the compactness of $\{\mathcal{G}_\mu\}$ in $C^{m, \kappa}(\T^N)$, this shows that the entire sequence converges uniformly to $\mathcal{G}$ as $\mu \rightarrow 0$. 

To show \eqref{eq:higher-block}, we first prove by induction that for $|\beta|\leq m$
\begin{equation}\label{eq:DzDbeta-X}
  \|D_Z\partial^\beta \G_\mu\|_q\le C_q|Z|^\alpha,
\end{equation}
where recall that $\tau_Z V(Y) = V(Y+Z)$ and $D_Z V(Y) = \tau_ZV(Y) - V(Y)$. The base case $|\beta| = 0$ has been shown in \eqref{eq:cell-phase-difference}. Assume \eqref{eq:DzDbeta-X} is true for $|\beta| \le m - 1$. Then for $|\beta| = m$, by applying $\partial^\beta$ to equation \eqref{eq:cell-equation}, we obtain formally
\begin{equation}\label{eq:Dbeta-X}
  (\mu^2 - \nabla_\lambda \cdot A \nabla_\lambda) \partial^\beta \X_\mu = \nabla_\lambda \cdot E_{\mu,\beta},
\end{equation}
where
\begin{equation*}
  E_{\mu,\beta} = \partial^\beta A+
  \sum_{0<\gamma\le\beta}\binom{\beta}{\gamma}
  (\partial^\gamma A)\partial^{\beta-\gamma} \G_\mu.
\end{equation*}
By induction, it is easy to show that
\begin{align}\label{eq:Gmu_Cm}
  \|\partial^\beta \mathcal{G}_\mu\|_q\leq C_q\quad \text{for }|\beta|\leq m,
\end{align}
using Lemma \ref{lem:torus-estimate} and the $C^{m,\alpha}$ regularity of $A$.
Further, applying $D_Z$ to equation \eqref{eq:Dbeta-X}, one can get
    \begin{equation*}
    (\mu^2 - \nabla_\lambda \cdot A \nabla_\lambda) 
    D_Z\partial^\beta \X_\mu
    =\nabla_\lambda \cdot\{ D_Z E_{\mu,\beta}
              +(D_Z A)\tau_Z \partial^\beta \G_\mu \}.
    \end{equation*}
    Using the $C^{m,\alpha}$ regularity of $A$ and the inductive assumption of $\G_\mu$, we have
    \begin{equation*}
        \| D_Z E_{\mu,\beta} \|_q + \| (D_Z A)\tau_Z \partial^\beta \G_\mu \|_q \le C|Z|^\alpha,
    \end{equation*}
    which, by Lemma \ref{lem:torus-estimate}, yields \eqref{eq:DzDbeta-X}. Moreover, in view of $\cH_\mu = A(I + \G_\mu)$, we obtain from \eqref{eq:DzDbeta-X} and \eqref{eq:Gmu_Cm} that for any $|\beta| \le m$
    \begin{align}\label{eq:Hmu_bound}
        \|\partial^\beta \cH_\mu\|_q \le C_q, \qquad \|D_Z\partial^\beta \cH_\mu\|_q
    \le C_q|Z|^\alpha.
    \end{align}
    
    Now, we prove \eqref{eq:higher-block} from \eqref{eq:DzDbeta-X} and \eqref{eq:Hmu_bound}. For simplicity, we focus on the estimate of $\G_\mu$, since the estimate of $\cH_\mu$ is similar. Note that the case $\ell = 0$ is trivial by \eqref{eq:phase-block-estimate}. For $\ell \ge 1$, recall that
    \begin{equation}\label{eq:DG-integral}
        \Delta_\ell \G_\mu(Y) = \int_{\R^N} \mathcal{K}_\ell(Z) \G_\mu(Y-Z) \dd Z,
    \end{equation}
    where the kernel $\mathcal{K}_\ell$ is given by \eqref{eq:kernel}. In view of \eqref{eq:phase-kernel-moments}, for any polynomial $P(Z)$, we have
    \begin{equation}\label{eq:moments-cancel}
        \int_{\R^N} \mathcal{K}_\ell(Z) P(Z) \dd Z = 0.
    \end{equation}
Since $\mathcal{G}_\mu\in C^m(\T^N)$, by the Taylor expansion at $Y$ we can write
    \begin{equation}\label{eq:TaylorExp}
    \begin{aligned}
        &\qquad\G_\mu(Y - Z) = \sum_{k=0}^{m} \frac{1}{k!} \nabla_Y^k \G_\mu(Y) : (-Z)^{\otimes k} \\
        & + \int_0^1 \int_0^{t_1} \cdots \int_0^{t_{m-1}} (\nabla^m \G_\mu(Y - t_m Z) - \nabla^m \G_\mu(Y)) : (-Z)^{\otimes m} \dd t_m \cdots \dd t_1.
    \end{aligned}
    \end{equation}
    Here for convenience, we have used the tensor notation (particularly $X^{\otimes k} = X\otimes X \dots \otimes X$), and $T:S$ means contraction between tensors $T$ and $S$. Note that the formula above holds for any $Z\in \R^N$. It follows from \eqref{eq:DG-integral}, \eqref{eq:moments-cancel} and \eqref{eq:TaylorExp} that
    \begin{equation*}
        \Delta_\ell \G_\mu(Y) = \int_{\R^N} \int_0^1 \int_0^{t_1} \cdots \int_0^{t_{m-1}} \mathcal{K}_\ell(Z) (D_{-t_m Z} \nabla^m\G_\mu(Y)):(-Z)^{\otimes m} \dd t_m \cdots \dd t_1 \dd Z.
    \end{equation*}
    Taking $L^q$ norm on both sides and using \eqref{eq:DzDbeta-X}, we have
    \begin{equation*}
    \begin{aligned}
        \| \Delta_\ell \G_\mu \|_q & \le \int_{\R^N} \int_0^1 \int_0^{t_1} \cdots \int_0^{t_{m-1}} |\mathcal{K}_\ell(Z)| t_m^\alpha |Z|^{m+\alpha} \dd t_m \cdots \dd t_1 \dd Z \\
        & \le C 2^{-(m+\alpha)\ell} \int_{\R^N} |\mathcal{K}(Z)||Z|^{m+\alpha} \dd Z \\
        & \le C 2^{-(m+\alpha)\ell}.
    \end{aligned}
  \end{equation*}
  The proof is complete.
\end{proof}

\begin{remark}
    The above lemma only shows that $\nabla_\lambda \X_\mu$ and $A(I + \nabla_\lambda \X_\mu)$ are uniformly bounded. However, $\X_\mu$ itself need not remain bounded as $\mu \to 0$, and thus we are not able to take limit of $\X_\mu$ directly.
\end{remark}
    
We point out that the endpoint estimate $q = \infty$ of \eqref{eq:G+H} can not be obtained here due to an essential loss from Bernstein's inequality. In fact, Lemma \ref{lem:corrector-existence} tells us that $\mathcal{G}$ belongs to the Besov space $B^{m+\alpha}_{q, \infty}(\T^N)$, which properly contains $B^{m+\alpha}_{\infty, \infty}(\T^N)=C^{m, \alpha}(\T^N)$ for $0<\alpha<1$. In order to recover the full exponent in the growth estimate of the true correctors, we need a more delicate low-frequency cutoﬀ estimate of physical orbits and a combinatorial counting argument.

We start with the equation of $\Delta_\ell\mathcal{G}$. Recall that
\begin{equation*}
    -\nlam \cdot (A (e_j + \G_j)) = 0 \qquad \text{in } \T^N.
\end{equation*}
Applying $\Delta_\ell$ to the equation, we get
\begin{equation}\label{eq:Delta-Gj}
    - \nlam \cdot (A \Delta_\ell \G_j ) = \nlam \cdot ( \Delta_\ell A e_j + [\Delta_\ell, A] \G_j ),
\end{equation}
where $[\Delta_\ell, A] G = \Delta_\ell (A G) - A \Delta_\ell G$ for $G \in L^2(\T^N; \R^d)$.

\begin{lemma}\label{lem:N}
    Let
    \begin{equation}\label{eq:def-Nell}
        \cN_\ell  := \Delta_\ell A + [\Delta_\ell, A] \G = \Delta_\ell \cH - A\Delta_\ell \G.
    \end{equation}
    Then
    \begin{equation}\label{eq:Nell-decay}
        \| \cN_\ell \|_\infty \le C 2^{-(m+\alpha) \ell}.
    \end{equation}
\end{lemma}

\begin{proof}
It is sufficient to prove the estimate for commutator $[\Delta_\ell, A] \G$, since the estimate for $\Delta_\ell A$ is simpler. The case $\ell = 0$ is trivial. For $\ell \ge 1$, note that
    \begin{equation*}
        [\Delta_\ell, A] \G(Y) = \int_{\R^N} \mathcal{K}_\ell(Z) (A(Y-Z) - A(Y)) \G(Y - Z) \dd Z.
    \end{equation*}
    Using an analogous formulation as \eqref{eq:TaylorExp}, we can write
    \begin{equation}
        A(Y-Z) - A(Y) = P_1(Z) + R_1(Z), \qquad \G(Y - Z) = P_2(Z) + R_2(Z),
    \end{equation}
    where $P_1, P_2$ are polynomials of degree at most $m$, such that, $$|P_1(Z)| \le C(|Z| + |Z|^m)\quad |P_2(Z)| \le C(1 + |Z|^m),$$
    and $$|R_1(Z)| \le C|Z|^{m+\alpha},\quad |R_2(Z)| \le C|Z|^{m+\kappa},\quad \kappa\in(0, \alpha).$$
    By \eqref{eq:moments-cancel}, we then have
    \begin{equation*}
    \begin{aligned}
        | [\Delta_\ell, A] \G(Y) |  & = \bigg| \int_{\R^N} \mathcal{K}_\ell(Z) (P_1 R_2 + P_2 R_1 + R_1 R_2) \dd Z \bigg| \\
        & \le \int_{\R^N} |\mathcal{K}_\ell(Z)| (|Z|^{1+m+\kappa} + |Z|^{m+\alpha} + |Z|^{2m+\alpha} + |Z|^{2m+\alpha+\kappa}) \dd Z \\
        & \le C 2^{-(m+\alpha) \ell}.
    \end{aligned}
    \end{equation*}
    This estimate is independent of $Y$, which ends the proof.
\end{proof}

\begin{lemma}\label{lem:potentials}
There are trigonometric polynomials $g_\ell=(g_\ell^i)$ and
$h_\ell=(h_\ell^{ijk})$, with $h_\ell^{ijk}=-h_\ell^{jik}$, such that
\begin{align}
 \nlam g_\ell&=\Delta_\ell \G,\label{eq:h-gradient}\\
 \sum_{i=1}^d\nabla_{\lambda,i}h_\ell^{ijk}
   &=(\Delta_\ell \cH)_{jk}-\Pi_\lambda((\Delta_\ell \cH)_{jk}),
   \label{eq:s-div}
\end{align}
and for every $2\leq q<\infty$
\begin{equation}\label{eq:orbit-est}
 \sup_Y\bigl[\Pi_\lambda(|\nlam g_\ell|^q)(Y)\bigr]^{1/q}
 +\sup_Y\bigl[\Pi_\lambda (|\nlam h_\ell|^q)(Y)\bigr]^{1/q}
 \le C_q2^{-(m+\alpha)\ell}.
\end{equation}
\end{lemma}
\begin{proof}
For $n\in\Z^N$ with $\xi_\lambda(n)\neq 0$, define
\begin{align}
 \widehat g_\ell^{\,i}(n)
 &=\frac{\xi_\lambda(n)\cdot\widehat{\Delta_\ell \G_i}(n)}
          {2\pi i|\xi_\lambda(n)|^2},\label{eq:h-formula}\\
 \widehat h_\ell^{\,ijk}(n)
 &=\frac{(\xi_\lambda(n))_i\widehat{(\Delta_\ell \cH_{jk})}(n)
        -(\xi_\lambda(n))_j\widehat{(\Delta_\ell \cH_{ik})}(n)}
         {2\pi i|\xi_\lambda(n)|^2}.
 \label{eq:s-formula}
\end{align}
If $\xi_\lambda(n)=0$, simply define $\widehat{g}_\ell(n)=\widehat{h}_\ell(n)=0$. It is easy to verify that $h_\ell^{ijk}$ is skew-symmetric in $i, j$ and \eqref{eq:h-gradient}--\eqref{eq:s-div} hold. 

By \eqref{eq:Delta-Gj} and \eqref{eq:h-gradient}, on each physical orbit $Y + \lambda \R^d$, $g^j_\ell(Y+\lambda x)$ satisfies
\begin{equation*}
 -\nabla_x\cdot\big(A(Y+\lambda x)
       \nabla_x[g^j_\ell(Y+\lambda x)]\big) 
   =\nabla_x\cdot \cN_{\ell,j}(Y+\lambda x).
\end{equation*}
Note that $g_\ell$ is a trigonometric polynomial. It follows from Lemma \ref{lem:physical-cz}, Lemma \ref{lem:N} as well as the Caccioppoli inequality that
\begin{align*}
  \|\nlam g^j_\ell (Y+\lambda\cdot)\|_{\underline L^q(B_T)}
  & \le C_q\{T^{-1}\|g_\ell\|_\infty + \| \cN_{\ell,j}(Y+\lambda\cdot)\|_{\underline L^q(B_{2T})} \} \\
  & \le C_q\{T^{-1}\|g_\ell\|_\infty + 2^{-(m+\alpha) \ell} \}.
\end{align*}
Letting $T \to \infty$ and using \eqref{eq:resonant-orbit-average}, we obtain
\begin{equation}\label{eq:hl-orbit-est}
    \sup_Y\bigl[\Pi_\lambda(|\nlam g_\ell|^q)(Y)\bigr]^{1/q} \le C 2^{-(m+\alpha) \ell}.
\end{equation}

On the other hand, $h_\ell^{ijk}$ satisfies
\begin{equation}\label{eq:sl-orbit}
 \Delta_\lambda h_\ell^{ijk}
 =\nabla_{\lambda,i}(\Delta_\ell \cH_{jk})
       -\nabla_{\lambda,j}(\Delta_\ell \cH_{ik}).
\end{equation}
Since we do not have a sharp $L^\infty$ bound of  $\Delta_\ell \cH$ as like $\mathcal{N}_\ell$, we instead use the identity \eqref{eq:def-Nell} and \eqref{eq:h-gradient} to get
\begin{align*}
    \|\Delta_\ell \cH (Y+\lambda\cdot)\|_{\underline L^q(B_T)} \le \| \cN_\ell \|_\infty + \| A \nlam g_\ell (Y+\lambda\cdot)\|_{\underline L^q(B_T)}.
\end{align*}
Thus, applying Lemma \ref{lem:physical-cz} to \eqref{eq:sl-orbit}, one can see
\begin{align*}
  \|\nlam h_\ell (Y+\lambda\cdot)\|_{\underline L^q(B_T)} & \le C \{ T^{-1} \| h_\ell \|_\infty + \|\Delta_\ell \cH (Y+\lambda\cdot)\|_{\underline L^q(B_{2T})} \} \\
                                                          & \le C \{ T^{-1} \| h_\ell \|_\infty + \| \cN_\ell \|_\infty + \| \nlam g_\ell (Y+\lambda\cdot)\|_{\underline L^q(B_{2T})} \}.
\end{align*}
Letting $T\to \infty$ and using \eqref{eq:Nell-decay} and \eqref{eq:hl-orbit-est}, we obtain
\begin{align*}
    \sup_Y\bigl[\Pi_\lambda (|\nlam h_\ell|^q)(Y)\bigr]^{1/q}
    \le C_q2^{-(m+\alpha)\ell}.
\end{align*}
\end{proof}

\begin{remark}\label{rmk:orbit-torus}
    The orbit estimate \eqref{eq:orbit-est} is stronger than the $L^q$ estimate \eqref{eq:G+H}. To see this, note that by \eqref{es:L1Linfty} 
    \begin{align*}
        \| \nabla_\lambda g_\ell \|_q^q = \int_{\T^N} \Pi_\lambda( |\nlam g_\ell|^q)(Y) \dd Y \le\sup_{Y\in \T^N} \Pi_\lambda(|\nlam g_\ell|^q)(Y).
    \end{align*}
Thus, \eqref{eq:G+H} can be deduced from \eqref{eq:orbit-est}.
\end{remark}

Now we are ready to construct correctors and flux correctors from $g_\ell$ and $h_\ell$. To construct correctors, a straightforward idea is to take the summation $\chi^Y(x) = \sum_{\ell \in \N} g_\ell(Y + \lambda x)$ restricted to each physical orbit $Y + \lambda \R^d$ with $Y \in \T^N$. By \eqref{eq:h-gradient}, this formally implies $\nabla_x \chi^Y(x) = \G(Y + \lambda x)$. However, since $g_\ell$ does not have a summable control, such a direct summation is not possible.

Our idea is to take advantage of decay estimates of the gradients. In fact, by Lemmas \ref{Bern}, \ref{lem:potentials} and Remark \ref{rmk:orbit-torus}, it holds that
\begin{equation}\label{eq:uniform-gradient-sum}
 \|\nlam g_\ell\|_\infty+\|\nlam h_\ell\|_\infty
 \le C_q2^{-(m+ \kappa)\ell}
\end{equation}
for any $\kappa \in (0,\alpha)$. This inspires us to define the true correctors and flux correctors by
\begin{align}
 \chi^Y(x)&:=\sum_{\ell\ge0}
        \{g_\ell(Y+\lambda x)-g_\ell(Y)\},\label{eq:def-chi}\\
 \phi^Y(x)&:=\sum_{\ell\ge0}
        \{h_\ell(Y+\lambda x)-h_\ell(Y)\}.\label{eq:def-phi}
\end{align}
Both series and their physical gradients converge uniformly on bounded domains, uniformly in $Y$, and
\begin{align}
  \nabla_x\chi^Y(x)=\sum_{\ell\geq 0}\nabla_\lambda g_\ell(Y+\lambda x),\quad
  \nabla_x\phi^Y(x)=\sum_{\ell\geq 0}\nabla_\lambda h_\ell(Y+\lambda x). \label{eq:gradchiphi}
\end{align}
In particular, these fields are $C^1$ in $x$, vanishing at zero, and satisfy
\begin{gather}
-\nabla_x \cdot [A_\lambda^Y(I+\nabla_x\chi^Y)] = 0,
    \label{eq:cell-exact} \\
    \Delta \phi^Y_{kij}=\partial_k[\cH_{ij}(Y+\lambda x)]-\partial_i[\cH_{kj}(Y+\lambda x)], \label{eq:fluxpossion}\\
    \sum_k\partial_k \phi^Y_{kij}
      =\cH_{ij}(Y+\lambda x)-(\overline A_\lambda^Y)_{ij},\label{eq:flux-exact}
    \end{gather}
    where $A_\lambda^Y(x)=A(Y+\lambda x)$ and
    \begin{align}\label{eq:def-Abar}
      \overline A_\lambda^Y:= \Pi_\lambda (\cH)(Y).
    \end{align}
By \eqref{eq:resonant-orbit-invariance}, $\overline A_\lambda^Y$ is constant on each orbit and therefore independent of $x$. If there are no nonzero integer resonances (i.e., $\xi_\lambda(n) \neq 0$ for any nonzero $n\in \Z^N$), it is also independent of $Y$. Besides, by a classical argument, one can show that 
\begin{align*}
    \zeta\cdot\overline A_\lambda^Y\zeta\geq \Lambda^{-1}|\zeta|^2~\text{ and }~ |\overline{A}_\lambda^Y\zeta|\leq \Lambda^3|\zeta|\qquad\text{for }\zeta\in \R^d. 
\end{align*}

Now we turn to the sublinear growth of $\chi^Y$ and $\phi^Y$ constructed in \eqref{eq:def-chi} and \eqref{eq:def-phi}. Note that
\begin{align}\label{eq:chi-g}
    \chi^Y(x) = x\cdot \int_0^1 \sum_{\ell \ge 0} \nlam g_\ell(Y + t\lambda x) \dd t.
\end{align}
The key is to control $\sum_{\ell \ge 0} \nlam g_\ell$. 

Define the smallest nonresonant frequency in each dyadic phase ball by
\begin{equation}\label{eq:nu}
    \nu_\ell:=\min\{|\lambda^Tn|:n\in\Z^N,\ |n|\le2^{\ell+1},\ \lambda^Tn \ne 0 \}.
\end{equation}
If the set is empty, we set $\nu_\ell=\infty$ and $\nu_\ell^{-1}=0$ for convention. 

\begin{lemma}\label{lem:improve-cancel}
Let the multiplier $S_k^\lambda$ be given by
\begin{align*}
  S_k^\lambda V(Y):=\sum_{n\in\Z^N}\vartheta(2^{-k}|\xi_\lambda(n)|)\widehat V(n)e^{2\pi i n\cdot Y}, \qquad k\in\Z,
\end{align*}
where $\vartheta$ is the cutoff given in Section \ref{sec.LP}. Then there exist a fixed integer $M$ and a constant $C$ such that
\begin{equation}\label{eq:lowpass}
 \|S_k^\lambda W\|_\infty\le C\sum_{\ell\ge0}2^{-(m+\alpha)\ell}\mathbf{1}_{\{\nu_\ell\le4M2^k\}},
 \qquad k\in\Z,
\end{equation}
for $W=\sum_{\ell\geq 0}\nabla_\lambda g_\ell$ or $\sum_{\ell\geq 0}\nabla_\lambda h_\ell$.
\end{lemma}
\begin{proof}
  Denote $W_\ell=\nabla_\lambda g_\ell$ or $\nabla_\lambda h_\ell$, respectively, and set $V_\ell=S_k^\lambda W_\ell$ for $k\in \Z$. Let $\widetilde{\mathcal{K}}^\vartheta$ be the inverse Fourier transform of $\vartheta(|\cdot|)$ in $\R^d$ and $\widetilde{\mathcal{K}}^\vartheta_k(x) = 2^{kd} \widetilde{\mathcal{K}}^\vartheta(2^k x)$. Then
\begin{align*}
    V_\ell(Y) = \int_{\R^d} \widetilde{\mathcal{K}}^\vartheta_k(x) W_\ell (Y - \lambda x) \dd x.
\end{align*}
By \eqref{eq:resonant-orbit-average} and Minkowski's inequality, we have
\begin{align*}
  \Pi_\lambda(|V_\ell|^q)(Y)\leq \|\widetilde{\mathcal{K}}^\vartheta\|_{L^1(\R^d)}^q\Pi_\lambda(|W_\ell|^q)(Y),
\end{align*}
which, together with \eqref{eq:orbit-est}, yields that 
\begin{equation}\label{eq:fl-orbit}
 \sup_Y[\Pi_\lambda (|V_\ell|^q)(Y)]^{1/q} \le C\sup_Y[\Pi_\lambda (|W_\ell|^q)(Y)]^{1/q} \le C_q2^{-(m+\alpha)\ell}.
\end{equation}

Choose an even integer $q=2M>\max\{2, N/(m+\alpha)\}$. Denote
\begin{equation*}
 L=\min\{\ell\ge0:\nu_\ell\le4M2^k\},
\end{equation*}
and, if this set is empty, set $L=\infty$. Note that $\nu_\ell$ decreases in $\ell$ and $\nu_{L-1}>4M2^k$, which yields $\one_{\{\nu_\ell\leq 4M 2^k\}}=\one_{\{\ell\geq L\}}$. If $0<|\lambda^Tn|\le4M2^k$, then $|\lambda^Tn|<\nu_{L-1}$ and, by \eqref{eq:nu}, $|n|>2^{L}$. Observing that the phase support of $V_\ell$ lies in $$I_{\ell,k} = \{n\in \Z^N: |n|\le2^{\ell+1}, 0<|\lambda^Tn|\le2^{k+1} \},$$ this also indicates that $V_\ell=0$ for $\ell<L$.

Now we deduce a pointwise estimate of $V_\ell$ from the $L^q$ orbit estimate \eqref{eq:fl-orbit}. Fix $Y$ and write
\begin{align*}
    V_\ell(Y + \lambda x) = \sum_{n\in I_{\ell, k}} a_n e^{2\pi i n \cdot (Y + \lambda x)} = \sum_{n\in I_{\ell,k}} a_n e^{2\pi in\cdot Y} e^{2\pi i \lambda^T n \cdot x}.
\end{align*}
Taking the $M$-th power, we have
\begin{equation}\label{eq:power-argument}
 [V_\ell(Y+\lambda x)]^M
       =\sum_{\xi\in \Xi} c_\xi(Y)e^{2\pi i\xi \cdot x},
\end{equation}
where the summation is taken over $\Xi=\{\lambda^T\tilde{n}: \tilde{n}\in \widetilde{I_{\ell, k}}\}$, and $\widetilde{I_{\ell, k}}=\{n_1+\cdots+n_M: n_j\in I_{\ell, k}\}$. For each $\tilde{n}\in\widetilde{I_{\ell, k}}$, it holds that $|\tilde{n}|\le M2^{\ell+1}$ and $|\lambda^T\tilde{n}|\le M2^{k+1}$. For different $\tilde{n}, \tilde{n}'\in \widetilde{I_{\ell, k}}$ such that $\lambda^T (\tilde{n} - \tilde{n}')\neq 0$, we know $0< |\lambda^T (\tilde{n} - \tilde{n}')|\le 4M2^k$ and further, by the discussion above, $|\tilde{n} - \tilde{n}'|> 2^L$. This tells us that $\tilde{n}\in \widetilde{I_{\ell, k}}$ with different physical frequencies are separated by $2^L$. In view of $|\tilde{n}|\leq M2^{\ell+1}$, we know the cardinal number of $\Xi$ is bounded by $C_{M,N} 2^{N(\ell - L)}$.

Applying the Cauchy-Schwarz inequality to \eqref{eq:power-argument} at $x=0$, we obtain
\begin{equation*}
 |V_\ell(Y)|^{2M}
 \le C_{M,N}2^{N(\ell-L)}\sum_{\xi\in\Xi}|c_\xi(Y)|^2
 =C_{M,N}2^{N(\ell-L)}\Pi_\lambda(|V_\ell|^{2M})(Y),
\end{equation*}
where in the last step we have used the fact
\begin{align*}
  \Pi_\lambda(|V_\ell|^{2M})(Y)&=\lim_{R\rightarrow\infty}\fint_{B_R}[V_\ell(Y+\lambda x)]^{2M} \dd x\\&=\lim_{R\rightarrow \infty}\fint_{B_R}\sum_{\xi,\xi'\in \Xi}c_\xi(Y)\overline{c_{\xi'}(Y)}e^{2\pi i(\xi-\xi')\cdot x}\dd x\\& = \sum_{\xi\in \Xi}|c_\xi(Y)|^2
\end{align*}
from \eqref{eq:resonant-orbit-average} and the Riemann-Lebesgue lemma. Taking supremum over $Y \in \T^N$ and using \eqref{eq:orbit-est}, we obtain
\begin{equation*}
 \|V_\ell\|_\infty
 \le C_q 2^{N(\ell-L)/q}2^{-(m+\alpha)\ell}.
\end{equation*}
Since $q>N/(m+\alpha)$, it follows that
\begin{equation*}
    \sum_{\ell\ge L} \|V_\ell\|_\infty \le
    2^{-(m+\alpha)L}\sum_{\ell\ge L}2^{-(m+\alpha-N/q)(\ell-L)}
    \le C2^{-(m+\alpha) L}.
\end{equation*}
Also recall that $V_\ell = 0$ for $\ell < L$. Thus,
\begin{align*}
    \| S_k^\lambda W\|_\infty\leq \sum_{\ell\geq L}\|V_\ell\|_\infty \le C 2^{-(m+\alpha) L}=C\sum_{\ell\geq 0}2^{-(m+\alpha)\ell}\one_{\{\ell\geq L\}}.
\end{align*}
This completes the proof.
\end{proof}

\begin{proof}[Proof of Theorem \ref{thm:correctors}]
  Let $\chi$ and $\phi$ be given by \eqref{eq:def-chi} and \eqref{eq:def-phi} with $Y=0$ respectively, which already satisfy \eqref{eq:corrector}--\eqref{eq:fluxcorrector} according to \eqref{eq:cell-exact}--\eqref{eq:fluxpossion}. Besides, $\chi(0)=0$, $\phi(0)=0$, and \eqref{eq:gradientchiphi} follows from \eqref{eq:uniform-gradient-sum} and \eqref{eq:gradchiphi}.


It remains to prove the sublinear growth \eqref{eq:corrector-growth}. We will focus on the estimate for $\chi$, as the estimate for $\phi$ is exactly the same. Let $k\in \Z$ be an integer to be determined. Observe that 
\begin{equation*}
  \mathrm{Id} = S_k^\lambda+\sum_{j>k}P_j^\lambda. 
\end{equation*}
We have from \eqref{eq:chi-g} that
\begin{align}\label{eq:chi-integral}
  \begin{split}
    \chi(x)= x\cdot \int_0^1 \mathcal{G} (t\lambda x) \dd t=x\cdot \int_0^1 \Big(S_k^\lambda \mathcal{G} (t\lambda x) + \sum_{j > k} P^\lambda_j \mathcal{G}(t\lambda x) \Big) \dd t.
  \end{split}
\end{align}
For $j > t$, set
\begin{align*}
    \widehat{T^\lambda_j \mathcal{G}}(n) = \frac{\rho(2^{-j} \xi_\lambda(n))  \xi_\lambda(n) \cdot \widehat{\mathcal{G}}(n)}{2\pi i |\xi_\lambda (n)|^2},
\end{align*}
which satisfies
\begin{equation}\label{eq:DTW}
    \nlam T^\lambda_j \mathcal{G} = P^\lambda_j \mathcal{G}.
\end{equation}
Note that $T_j^\lambda$ is well-defined as $\rho(2^{-j} \xi_\lambda(n))$ is supported in $\{ 2^{j-1} < |\xi_\lambda(n)| \le 2^{j+1} \}$. Setting 
\begin{align*}
    \widetilde{\mathcal{K}} =\mathcal{F}^{-1}\Big[ \frac{\rho(\varsigma)  \varsigma }{2\pi i |\varsigma|^2}\Big]
\end{align*}
via the inverse Fourier transform in $\R^d$, we have as in \eqref{iden:convolution}
\begin{align*}
    T_j^\lambda \mathcal{G}(Y)=2^{-j}\int_{\R^d}\widetilde{\mathcal{K}}(y)\mathcal{G}(Y-\lambda y/2^j)\dd y.
\end{align*}
Besides, since $\vartheta(2^{-j-1} |\xi_\lambda(n)|) = 1$ on the support of $\rho(2^{-j}\xi_\lambda(n))$, we have $T_j^\lambda \mathcal{G}=T_j^\lambda S_{j+1}^\lambda \mathcal{G}$. Consequently,
\begin{align*}
    T^\lambda_j \mathcal{G}(Y) =2^{-j}\int_{\R^d}\widetilde{\mathcal{K}}(y)S_{j+1}^\lambda\mathcal{G}(Y-\lambda y/2^j)\dd y,
\end{align*}
which yields
\begin{equation}\label{eq:Tj-infty}
    \| T^\lambda_j \mathcal{G}\|_\infty \le C2^{-j} \| S_{j+1}^\lambda \mathcal{G} \|_\infty.
\end{equation}

Now, taking \eqref{eq:DTW} into \eqref{eq:chi-integral}, we have
\begin{align*}
    \chi(x) = x\cdot \int_0^1 S_k^\lambda \mathcal{G} (t\lambda x)  \dd t + \sum_{j > k} [T^\lambda_j \mathcal{G}(\lambda x) - T^\lambda_j \mathcal{G}(0)].
\end{align*}
Recall that $\mathcal{G}=\sum_{\ell\geq 0}\nabla_\lambda g_\ell$. Then, Lemma \ref{lem:improve-cancel}, together with \eqref{eq:Tj-infty}, implies
\begin{align*}
    \sup_{|x| \le R} |\chi(x)| \le C R \sum_{\ell \ge 0} 2^{-s\ell} \mathbf1_{\{\nu_\ell\le4M2^k\}}
   +C \sum_{j>k}2^{-j} \sum_{\ell \ge 0} 2^{-s\ell} \mathbf1_{\{\nu_\ell\le4M2^{j+1}\}}.
\end{align*}
Choose $k \in \Z$ such that $2^k \le R^{-1} < 2^{k+1}$. By direct calculation, for every $\ell \ge 0$,
\begin{equation*}
 R\mathbf1_{\{\nu_\ell\le4M2^k\}}
   +\sum_{j>k}2^{-j}\mathbf1_{\{\nu_\ell\le4M2^{j+1}\}}
 \le C_M\min\{R,\nu_\ell^{-1}\}.
\end{equation*}
By the definitions of $\cJ^\lambda_\ell$ and $\nu_\ell$, if $\cJ_\ell^\lambda$ is nonempty, then
$2^{\min\cJ_\ell^\lambda-1}\leq \nu_\ell\leq 2^{\min\cJ_\ell^\lambda+1}$. Therefore, in this case
\begin{equation*}
 \min\{R,\nu_\ell^{-1}\}
 \le C\sum_{j\in\cJ_\ell^\lambda}\min\{R,2^{-j}\},
\end{equation*}
and further
\begin{equation}\label{eq:sublinear-rate0}
\begin{aligned}
    \sup_{|x| \le R} |\chi(x)| \le  C \sum_{\ell \ge 0} 2^{-s\ell} \min\{R,\nu_\ell^{-1}\} \le C \sum_{\ell \ge 0} 2^{-s\ell} \sum_{j\in\cJ_\ell^\lambda}\min\{R,2^{-j}\}.
\end{aligned}
\end{equation}
Note that if $\mathcal{J}_\ell^\lambda$ is empty, then $\nu_\ell=\infty$ and the estimate \eqref{eq:sublinear-rate0} is also valid.

    Observe that
    \begin{equation*}
        \Psi(t):= \int_t^\infty \psi(r) \frac{\dd r}{r} = 
        \begin{cases}
            2 -t, \quad & 0<t\le 1, \\
            t^{-1}, \quad & t > 1.
        \end{cases}
    \end{equation*}
    Hence, $\Psi(t) \simeq \min\{ 1, t^{-1} \}$ and
    \begin{equation*}
        \int_R^\infty \psi(2^j r) \frac{\dd r}{r} = \Psi(2^j R) \simeq \min \{ 1, (2^j R)^{-1} \}.
    \end{equation*}
    It follows that
    \begin{equation}\label{eq:equi-rates}
        \begin{aligned}
            \sum_{\ell \ge 0} 2^{-s\ell} \sum_{j\in\cJ_\ell^\lambda}\min\{R,2^{-j}\} & = R \sum_{\ell \ge 0} 2^{-s \ell} \sum_{j\in \cJ_\ell^\lambda}  \min \{1, (2^j R)^{-1} \} \\
            & \simeq CR \int_R^\infty \sum_{\ell \ge 0} 2^{-s \ell} \sum_{j\in \cJ_\ell^\lambda} \psi(2^j r) \frac{\dd r}{r} \\
            & = CR \int_R^\infty \omega_{\lambda,s}(r) \frac{\dd r}{r}.
        \end{aligned}
    \end{equation}
    Together with \eqref{eq:sublinear-rate0}, we get
    \begin{equation*}
        \sup_{Y\in \T^N} \sup_{|x|\le R} |\chi^Y(x)| \le CR \int_R^\infty \omega_{\lambda,s}(r) \frac{\dd r}{r}.
    \end{equation*}
    This proves \eqref{eq:corrector-growth}.

    Finally, one can also show the correctors and flux correctors satisfying \eqref{eq:corrector} and \eqref{eq:fluxcorrector} with the sublinear growth \eqref{eq:corrector-growth} are unique under the normalization at zero. This follows from the zeroth-order Liouville theorem (see Corollary \ref{cor:zero-liouville}).
\end{proof}

\begin{remark}
    We emphasize that Lemma \ref{lem:improve-cancel} is the key ingredient to obtain the optimal growth rate without any loss. If a small loss is allowed, we can actually give a much simpler proof. In fact, by turning \eqref{eq:h-formula} into a convolution integral, we can show
    \begin{equation*}
        \| g_\ell \|_\infty \le C \nu_\ell^{-1} \| \Delta_\ell \G \|_\infty \le C\nu_\ell^{-1} 2^{-(m+\kappa)\ell},
    \end{equation*}
    for any $\kappa \in (0,\alpha)$. Combined with \eqref{eq:uniform-gradient-sum}, we have
    \begin{equation*}
        |g_\ell(Y + \lambda x) - g_\ell(Y)| \le C \min\{ \nu_\ell^{-1},  |x| \} 2^{-(m+\kappa) \ell}.
    \end{equation*}
    Summing over $\ell$ leads to the growth rate for $\chi^Y$ with modulus $\omega_{\lambda, m+\kappa}$, which is suboptimal.
\end{remark}

As a straightforward consequence of the sublinear growth of the correctors, one can deduce the first-order Liouville theorem and large-scale $C^{1,1-}$ regularity as in \cite{GNO20}. The results in \cite{GNO20} were established under stochastic homogenization setting; but the argument there identically applies to the deterministic setting, using the notion of minimal radius. We emphasize that though the argument itself is quantitative, only a qualitative sublinear growth rate is needed. The following corollary is parallel to \cite[Corollary 1]{GNO20}. The proof will be omitted.

\begin{corollary}[First-order Liouville theorem]\label{thm:s6-Liouville}
Let $u$ be an entire $A_\lambda$-harmonic function in $\R^d$. If, for some
$0\le\gamma<1$ and $C_u<\infty$,
\begin{equation}\label{eq:s6-u-growth}
 \anorm{u}{B_R}\le C_uR^{1+\gamma},\qquad R\ge1,
\end{equation}
then there exist $a \in\R$ and $p\in\R^d$ such that
\begin{equation}\label{eq:s6-Liouville}
 u(x) = a + p\cdot(x + \chi(x)).
\end{equation}
Thus, this space of entire solutions has dimension $d+1$.
\end{corollary}

\subsection{Diophantine conditions and bounded correctors}
\label{subsec:s6-bounded}

\begin{theorem}[Diophantine growth rates]\label{thm:diophantine}
Suppose that $A$ satisfies the same assumptions as in Theorem \ref{thm:correctors} and $\lambda$ satisfies the Diophantine condition \eqref{eq:diophantine} with $\beta>0$. Then for $R\ge2$,
\begin{equation}\label{eq:Diophantine-rates}
 \sup_{|x|\le R}(|\chi(x)|+|\phi(x)|)
 \le C
 \begin{cases}
   1,&\text{if }\beta<s,\\
   \log_2 R,&\text{if }\beta=s,\\
   R^{1-s/\beta},&\text{if }\beta>s,
 \end{cases}
\end{equation}
where $s = m+\alpha$ and the constant $C$ may also depend on $c,\beta$ in \eqref{eq:diophantine}.
\end{theorem}

\begin{proof}
We focus on the estimate of $\chi$, since the other is the same. Recall the sublinear growth rate given by \eqref{eq:sublinear-rate0}
\begin{equation}\label{eq:nearest-growth}
    \sup_{|x| \le R} |\chi(x)|  \le  C \sum_{\ell \ge 0} 2^{-s\ell} \min\{R,\nu_\ell^{-1}\},
\end{equation}
where $\nu_\ell$ is defined by \eqref{eq:nu}. The Diophantine condition \eqref{eq:diophantine} gives $\nu_\ell\ge c_0 2^{-\beta\ell}$, where $c_0=c2^{-\beta}$.
Splitting the sum in \eqref{eq:nearest-growth} at
$\ell = L$, we get
\begin{equation*}
\begin{aligned}
    \sum_{\ell\ge0}2^{-s\ell}\min\{R,\nu_\ell^{-1}\}
    & \le c_0^{-1}\sum_{\ell=0}^L2^{(\beta-s)\ell} +R\sum_{\ell>L}2^{-s\ell} \\
    & \le c_0^{-1}\sum_{\ell=0}^L2^{(\beta-s)\ell} + \frac{ R 2^{-sL}}{2^s-1}.
\end{aligned}
\end{equation*}
Taking $L=\max\{1, \lfloor\beta^{-1} \log_2 R\rfloor\}$ and considering the three cases according to the relative magnitudes of $\beta$ and s, we obtain \eqref{eq:Diophantine-rates}. 
\end{proof}

\begin{remark}
    Under the Diophantine condition \eqref{eq:diophantine}, the existence of bounded correctors was known; see \cite{Koz78} and a more recent generalization \cite{AGK16} to almost-periodic homogenization. However, these results require unspecified higher-order regularity of coefficients. A careful examination of the argument in \cite{AGK16} shows that $A$ needs $C^m$ regularity with $m > N(\beta  +1)$; while the argument in \cite{Koz78} needs $m > \beta + N/2$. Theorem \ref{thm:diophantine} above shows that the optimal regularity assumption for coefficients is $m + \alpha > \beta$, independent of the dimension $N$.
\end{remark}

\subsection{Optimality}
In this section, we provide concrete examples showing that the sublinear growth rates established in Theorem \ref{thm:correctors} are generally optimal and can be exactly attained in some typical situations.

\begin{proposition}
\label{prop:optimal-example}
Let $m\ge0$, $0<\alpha<1$, $s=m+\alpha$, and let
$\beta >s$ be an integer. Set $N=\beta+2$, $\vartheta=2^{1/(\beta + 1)}$, and fix
\begin{equation}\label{eq:s6-example-lambda}
 \lambda=
 \begin{pmatrix}
 1&0\\ \vartheta&0\\ \vdots&\vdots\\ \vartheta^\beta&0\\0&1
 \end{pmatrix}\in\R^{N\times2}.
\end{equation}
There is a uniformly elliptic, symmetric $A\in C^{m,\alpha}(\T^N;\R^{2\times2})$ such that for all sufficiently large $R$,
\begin{equation}\label{eq:example-sharp}
     \sup_{|x|\le R}| \chi(x)|
    \simeq \sup_{|x|\le R}|\phi(x)|
    \simeq R \int_R^\infty\omega_{\lambda,s}(r)\,\frac{\dd r}{r}
    \simeq R^{1-s/\beta}.
\end{equation}
\end{proposition}
\begin{proof}
Let $n' = (n_1,\dots, n_{\beta+1})$. Define $P(z) = \sum_{j=0}^\beta n_{j+1} z^j$. Let $Q(z) = z^{\beta+1} - 2$. Thus $\vartheta$ is a simple root of $Q(z)$. By Eisenstein's criterion, $Q(z)$ is irreducible, which implies $P(\vartheta) \neq 0$ for any $n' \neq 0$. Hence, by Liouville's inequality \cite[Theorem A.1]{B04}, we have
\begin{equation}
    |P(\vartheta)| \ge C_\beta \Big(\sup_{0\le j \le \beta} |n_{j+1}| \Big)^{-\beta} \ge C_\beta |n'|^{-\beta}.
\end{equation}
Thus, for any $n = (n', n_{\beta+2}) \in \Z^{N} \setminus \{ 0 \}$, 
\begin{equation}\label{eq:algebraic-lowerbd}
|\lambda^T n| = |(P(\vartheta), n_{\beta+2})| \ge C_\beta |n|^{-\beta}.
\end{equation}

For the reverse estimate, let $M \ge 1$ and take the fractional parts of
\begin{equation*}
 \sum_{j=1}^{\beta} n_{j+1} \vartheta^j,
 \qquad (n_2,\ldots,n_{\beta+1})\in\{0,\ldots,M\}^{\beta}.
\end{equation*}
There are $(M+1)^\beta$ fractional parts and only $M^\beta$ intervals of length $M^{-\beta}$ in $[0,1)$. By pigeonhole principle, there are two fractional parts in the same interval of length $M^{-\beta}$. This means that there are $(n_2,\dots, n_{\beta+1})$ and $(\widetilde{n}_2,\dots, \widetilde{n}_{\beta+1})$ such that
\begin{equation*}
    \sum_{j=1}^{\beta} (n_{j+1} - \widetilde{n}_{j+1}) \vartheta^j = \text{ some integer } s_0 + O(M^{-\beta}).
\end{equation*}
Taking $n^* = (-s_0, n_{2} - \widetilde{n}_{2}, \dots, n_{\beta+1} - \widetilde{n}_{\beta+1}, 0) \in \Z^{N}$, we then have
\begin{equation}
    |\lambda^T n^*| \le CM^{-\beta}.
\end{equation}
Clearly, $|n^*| \le C_\beta M$. Hence, as $M$ varies, we obtain together with \eqref{eq:algebraic-lowerbd}
\begin{equation}\label{eq:algebraic-min}
 \min_{0<|n|\le M}|\lambda^Tn|
       \simeq_\beta M^{-\beta}, \qquad M \ge 1.
\end{equation}

In particular, by taking $M = 2^\ell$ with $\ell \ge 0$,
\begin{equation}\label{eq:nu-bound}
 \nu_\ell \simeq \min_{j\in \cJ_\ell^\lambda} 2^j \simeq 2^{-\beta \ell},\qquad
 \cJ_\ell^{\lambda}
       \subset[-\beta \ell-C_\beta,\ell+C_\beta]\cap\Z,
\end{equation}
and there is one $j \in \cJ_\ell^{\lambda}$  satisfying $2^j\simeq 2^{-\beta \ell}$.

Now, we show, for every $\eta \in (0, \beta)$ (including $\eta = s = m+\alpha$),
\begin{equation}\label{eq:example-tail}
    R\int_R^\infty\omega_{\lambda,\eta}(r)\,\frac{\dd r}{r}
       \simeq_{\beta,\eta}R^{1-\eta/\beta},
       \qquad R\ge2.
\end{equation}
We use the equivalent formula \eqref{eq:equi-rates} to prove this.
For a given $R \ge 2$, choose an integer $L\ge0$ with
$2^{\beta L}\le R<2^{\beta (L+1)}$. By \eqref{eq:nu-bound} and geometric summation, we have
\begin{equation*}
    \sum_{j\in\cJ_\ell^{\lambda}}\min\{R,2^{-j}\}
    \le C_\beta
    \begin{cases}
    2^{\beta \ell},&\ell\le L,\\
    R(1+\ell-L),&\ell>L,
    \end{cases}
\end{equation*}
and
\begin{equation*}
    \sum_{j\in\cJ_\ell^{\lambda}}\min\{R,2^{-j}\}
    \ge c_\beta
    \begin{cases}
    2^{\beta \ell},&\ell\le L,\\
    R,&\ell>L.
    \end{cases}
\end{equation*}
Consequently,
\begin{equation}
\begin{aligned}
    \sum_{\ell \ge 0}^\infty 2^{-\eta \ell } \sum_{j\in\cJ_\ell^{\lambda}}\min\{R,2^{-j}\} & \le C\sum_{0\le \ell \le L} 2^{(\beta - \eta) \ell} + \sum_{\ell > L} 2^{-\eta \ell} R(1+\ell - L) \\
    & \le C(2^{(\beta-\eta) L} + R 2^{-\eta L})  \le CR^{1-\eta/\beta},
\end{aligned}
\end{equation}
and
\begin{equation}
\begin{aligned}
    \sum_{\ell \ge 0}^\infty 2^{-\eta \ell } \sum_{j\in\cJ_\ell^{\lambda}}\min\{R,2^{-j}\} & \ge c\sum_{0\le \ell \le L} 2^{(\beta - \eta) \ell} + c\sum_{\ell > L} 2^{-\eta \ell} R \\
    & \ge c(2^{(\beta-\eta) L} + R 2^{-\eta L})  \ge c R^{1-\eta/\beta}.
\end{aligned}
\end{equation}
These estimates combined lead to \eqref{eq:example-tail}.

Next, we will construct a $C^{m,\alpha}$ coefficient matrix $A_\lambda$ whose corresponding correctors and flux correctors grow exactly as $R^{1-s/\beta}$, where $s = m + \alpha$.
Fix a sufficiently large integer $Q$. For $k \ge 1$, applying \eqref{eq:algebraic-min} to $M = Q^k$, and
changing signs if necessary, we see that there exists $n_k = (n_{k,j}) \in\Z^{\beta+1}$ such that
\begin{equation}\label{eq:fixed-modes}
 cQ^k\le |n_k|\le CQ^k,\qquad
 \delta_k:= \sum_{j=0}^\beta n_{k,j+1} \vartheta^j >0,\qquad
 cQ^{-\beta k}\le\delta_k\le CQ^{-\beta k}.
\end{equation}
The lower bound for $|n_k|$ follows from the algebraic lower
bound \eqref{eq:algebraic-lowerbd} and the upper bound for $\delta_k$.
All constants here depend only on $\beta$.
Increase $Q$ so that $\delta_{k+1}<\delta_k$ and, with
$B_k=Q^{-sk}/\delta_k$,
\begin{equation}\label{eq:s6-amplitude-separation}
 B_{k-1}\le B_k/3,\qquad k\ge2.
\end{equation}
This is possible because $\beta >s$ ($Q^{\beta -s} \ge 3C/c$ is sufficient).

Writing $Y'=(Y_1,\ldots,Y_{\beta+1})$ and $Y = (Y', Y_N)$, define
\begin{equation}\label{eq:s6-example-fg}
 f(Y) 
=f(Y')
=\varepsilon\sum_{k\ge1}Q^{-sk}
                     \cos(2\pi n_k\cdot Y'),
\end{equation}
where $\varepsilon>0$ is small and fixed, and $n_k$ are given in \eqref{eq:fixed-modes}.
All derivative series of $f$ through order $m$ converge uniformly.
For $|\gamma|=m$ and $0<h=|Z|\le1$, using the size of $n_k$ in \eqref{eq:fixed-modes},
\begin{equation*}
 |\partial^\gamma f(Y'+Z)-\partial^\gamma f(Y')|
 \le C\varepsilon\sum_{k\ge1}
          Q^{-\alpha k}\min\{1,Q^kh\}
 \le C\varepsilon h^\alpha.
\end{equation*}
The last inequality follows by splitting the sum at $Q^kh \simeq 1$.
Thus $f\in C^{m,\alpha}(\T^N)$ and $\|f\|_\infty \leq \frac{1}{2}$
after decreasing $\varepsilon$.
Set
\begin{equation}\label{eq:s6-example-A}
 A(Y)=\operatorname{diag}
       \left(\frac{1}{1+f(Y')},
             1+f(Y')\right), \qquad Y = (Y', Y_N) \in \T^N.
\end{equation}
This $2\times 2$ matrix is symmetric (diagonal), uniformly elliptic, and $C^{m,\alpha}(\T^N; \R^{2\times 2})$. This construction of coefficient matrix is known; see e.g., \cite[Section 8.7]{Mil02}.
The diagonal ratios are chosen so that the candidate corrector gradient cancels each denominator, leaving a divergence-free flux.
More precisely, for the coefficient matrix $A_\lambda(x) = A(\lambda x)$, the exact gradient and flux fields are
\begin{equation*}
    \G=\operatorname{diag}(f,0),\qquad
    \cH = A(I + \G )=\operatorname{diag}(1,1+f),\qquad
    \overline A_{\lambda}=I.
\end{equation*}

The normalized (at zero) correctors are
\begin{equation}\label{eq:example-chi}
 \chi_1(x)=\frac{\varepsilon}{2\pi}
       \sum_{k\ge1}\frac{Q^{-sk}}{\delta_k} \sin(2\pi\delta_kx_1),
       \qquad \chi_2(x) = 0.
\end{equation}
Clearly, $\nabla \chi_1(x) = \G_1(\lambda x) = (f(\lambda x), 0)$ and $\nabla \chi_2(x) = \G_2(\lambda x) = (0, 0)$.

The nontrivial components of flux correctors are
\begin{equation}\label{eq:example-fluxcorrector}
 \phi_{211}=\chi_2,\quad
 \phi_{121}=-\chi_2,\quad
 \phi_{122}=\chi_1,\quad
 \phi_{212}=-\chi_1.
\end{equation}
These formulas are derived from equation \eqref{eq:fluxpossion} and the uniqueness.

Now we verify the growth rates. For the corrector $\chi_1$, by the size of $\delta_k$ in \eqref{eq:fixed-modes}, the upper bound is given by
\begin{equation*}
 \sup_{|x|\le R}|\chi_1(x)|
 \le C\varepsilon\sum_{k\ge1}Q^{-sk}\min\{R,Q^{\beta k}\}
 \le CR^{1-s/\beta}.
\end{equation*}
For the lower bound set $R_k=(4\delta_k)^{-1}$.
At $x=(R_k,0)$ the $k$-th sine is one, and all later sines
are positive. Equation \eqref{eq:s6-amplitude-separation} gives
$\sum_{j<k}B_j\le B_k/2$, hence
\begin{equation*}
 \chi_1(R_k,0) \ge \frac{\varepsilon}{2\pi} (B_k - \sum_{\ell = 1}^{k-1} B_\ell ) \ge \frac{\eps}{4\pi} B_k \ge c R_k^{1-s/\beta}.
\end{equation*}
Moreover, since $R_k\simeq Q^{\beta k}$ and
$R_{k+1}/R_k\le C Q^\beta$,
the monotonicity of the supremum in $R$ therefore extends the lower bound
to every sufficiently large $R$. This, as well as  \eqref{eq:example-fluxcorrector}, implies
\begin{equation}
    \sup_{|x| \le R} |\chi(x)| \simeq \sup_{|x| \le R} |\phi(x)| \simeq R^{1-s/\beta}.
\end{equation}
Together with \eqref{eq:example-tail}, this completes the proof.
\end{proof}

\begin{remark}
    In the above example, we only consider integer $\beta > m + \alpha$ with $\alpha \in (0,1)$. For $\beta \le m$, Theorem \ref{thm:diophantine} directly implies the boundedness of correctors and flux correctors. Hence, the above example also shows that Theorem \ref{thm:diophantine} is generally optimal in terms of both the Diophantine condition and regularity of coefficients (we do not claim the sharpness for the endpoint cases that $\alpha = 1$ or $\beta = m+\alpha$, for which logarithmic improvement is possible).
\end{remark}

\bibliography{ref}
\bibliographystyle{abbrv}

\end{document}